\documentclass[10pt]{article}
\usepackage{amssymb}
\usepackage{amsmath}
\usepackage[dvips]{graphics}
\usepackage{epsfig}

\begin{document}
\title{Periodic Jacobi Operators with Complex Coefficients}

\author{
\textbf{Vassilis G. Papanicolaou}
\\\\
Department of Mathematics
\\
National Technical University of Athens
\\
Zografou Campus, 157 80, Athens, GREECE
\\
\underline{\tt papanico@math.ntua.gr}}
\maketitle

\begin{abstract}
We present certain results on the direct and inverse spectral theory of the Jacobi operator with complex periodic coefficients. For instance,
we show that any $N$-th degree polynomial whose leading coefficient is $(-1)^N$ is the Hill discriminant of finitely many discrete $N$-periodic
Schr\"{o}dinger operators (Theorem 1). Also, in the case where the spectrum is a closed interval we prove a result (Theorem 5) which is the analog
of Borg's Theorem for the non-self-adjoint Jacobi case.
\end{abstract}

\\


\section{Introduction}
We consider the periodic Jacobi (or discrete Hill-type operator) $L$ defined as
\begin{equation}
(L w)(n) := a(n) \, w(n+1) + a(n-1) \, w(n-1) + b(n)\, w(n),
\qquad
n \in \mathbb{Z},
\label{A3}
\end{equation}
where the coefficients  $a(n)$ and $b(n)$ are complex-valued and periodic functions of period $N \in \mathbb{N} := \{1, 2, \ldots\}$ with
\begin{equation}
a(n) \ne 0
\qquad \text{for all }\;
n \in \mathbb{Z}.
\label{A3a}
\end{equation}
Let us mention that if $a(n_0) = 0$ for some $n_0$, then $a(n_0 + kN) = 0$ for all $k \in \mathbb{Z}$ and $L$ splits as
$L = \bigoplus_{k \in \mathbb{Z}} A$, where
$A$ is a linear operator acting on an $N$-dimensional space, which can be considered as a degenerate case (e.g., the spectrum of $L$ consists of
at most $N$ eigenvalues of infinite multiplicity).

In the special case where $a(n) \equiv -1$ the operator $L$ becomes the one-dimensional discrete periodic Schr\"{o}dinger (or discrete Hill)
operator with potential $b(n)$. Unlike the continuous case, here there is no Liouville-type transformation which transforms the general operator $L$
of \eqref{A3} to a discrete periodic Schr\"{o}dinger operator (in fact, even in the continuous case, in the presence of complex coefficients the
Liouville transformation becomes problematic).

From now on, unless otherwise stated, without loss of generality we will normalize $a(n)$ so that
\begin{equation}
\prod_{j=1}^N a(j) = (-1)^N
\label{DA0}
\end{equation}
(starting with an arbitrary $a(n) \ne 0$ we can always do this normalization by replacing $L$ by $cL$, where $c$ is such that
$\prod_{j=1}^N c a(j) = c^N \prod_{j=1}^N a(j) = (-1)^N$).

Since $a(n)$ and $b(n)$ are $N$-periodic, they can be expanded as
\begin{equation}
a(n) = A_0 + \sum_{k=1}^{N-1} A_k \omega_N^{kn}
\quad \text{and} \quad
b(n) = B_0 + \sum_{k=1}^{N-1} B_k \omega_N^{kn},
\qquad
\omega_N := e^{2\pi i/N},
\label{DA1}
\end{equation}
where $A_0, A_1, \ldots, A_{N-1}, B_0, B_1, \ldots, B_{N-1} \in \mathbb{C}$ (this is a Fourier-style expansion). It is not hard to check the
orthogonality relation
\begin{equation}
\sum_{n=0}^{N-1} \omega_N^{jn} \, \bar{\omega}_N^{kn} = \sum_{n=0}^{N-1} e^{2(j-k)n\pi i/N} = N \delta_{jk}
\qquad \text{for }\; j,k = 0, 1, \ldots, N-1,
\label{DA2}
\end{equation}
where the bar denotes complex conjugation (thus $\bar{\omega}_N = \omega_N^{-1}$) and $\delta_{jk}$ is the Kronecker delta.
Using \eqref{DA2} in \eqref{DA1} yields
\begin{equation}
A_k = \frac{1}{N} \sum_{n=0}^{N-1} a(n) \, \bar{\omega}_N^{kn}
\qquad \text{and} \qquad
B_k = \frac{1}{N} \sum_{n=0}^{N-1} b(n) \, \bar{\omega}_N^{kn},
\label{DA3a}
\end{equation}
in particular (for $k=0$),
\begin{equation}
\sum_{n=0}^{N-1} a(n) = A_0 N
\qquad \text{and} \qquad
\sum_{n=0}^{N-1} b(n) = B_0 N.
\label{DA3}
\end{equation}

The continuous analog of $L$ is the operator $(H w)(x) := -w''(x) + V(x) w(x)$, where $V(x + b) = V(x)$. In the case where $V(x)$ is real-valued
there is a huge number of papers devoted to the spectral theory of $H$. But, even for the case where $V(x)$ is not real and, consequently, $H$ is not self-adjoint, there is an extensive amount of literature (see, e.g., \cite{G}, \cite{G-H}, \cite{G-T}, \cite{G-W}, \cite{G-U1}, \cite{G-U2},
\cite{P-T1}, \cite{P-T2}, \cite{P-T3}, \cite{RB}, \cite{S-T}, \cite{S}, \cite{T1}, \cite{T2}, \cite{T3}, \cite{T4}, \cite{V-TD}, as well as the
references therein). As for the discrete periodic case, there is, too, a considerable amount of
literature (e.g., \cite{D-T}, \cite{G-H-T}, \cite{Ho}, \cite{L}, \cite{N1}, \cite{N2}, \cite{P}, \cite{T}, \cite{Z}, and the references therein, as well as Barry Simon's encyclopedia \cite{Si}).

The recent emergence of the $\mathcal{PT}$-Symmetric Quantum Theory (see, e.g., \cite{B}) provides a strong motivation for studying non-self-adjoint
Schr\"{o}dinger-type operators (``non-Hermitian Hamiltonians" in the physicists' terminology), especially in the case where their spectra are real.

Let us recall that in the continuous case, if $V(x) \in L_{\text{loc}}^2(\mathbb{R})$ is real-valued, a famous
theorem of Borg \cite{Bo} states that $\sigma(H) = [\,0, \infty)$ if and only if $V(x) = 0$ a.e.\,. However, in the case of a nonreal $V(x)$ the
situation is very different since, as it has been shown by Gasymov \cite{G} (see also \cite{G-U2}), if
\begin{equation}
V(x) = \sum_{k=1}^{\infty} c_k e^{ikx},
\qquad \text{with }\;
\sum_{k=1}^{\infty} |c_k| < \infty,
\label{AA0}
\end{equation}
then the Hill discriminant of $H$ is $2\cos(2\pi\sqrt{\lambda})$ and, consequently, $\sigma(H) = [\,0, \infty)$.
Actually Gasymov's result can be easily explained, at least formally, by applying the substitution $z := e^{ix}$ to
the equation $(H w)(x) = \lambda w(x)$ and then invoking the standard Fuchsian theory of linear ordinary differential equations.
An interesting question here is whether, under some smoothness requirements, any potential $V(x)$ whose spectrum is $[\,0, \infty)$ must be a
``Gasymov potential," i.e. of the form given by \eqref{AA0}, or the complex conjugate of a Gasymov potential.

Clearly our discrete operator $L$ is bounded on $l^2(\mathbb{Z})$, and hence the $l^2(\mathbb{Z})$-spectrum $\sigma(L)$ of $L$ is a compact subset of
$\mathbb{C}$. In the present article, inspired by the aforementioned remarks on the continuous case, we examine, among other things, the somehow
simplest case regarding the spectrum, namely the case where $\sigma(L)$ is a closed interval.

\section{Review of the spectral theory of $L$}
Naturally, the spectral theory of $L$, acting on $l^2(\mathbb{Z})$, is studied via the equation
\begin{equation}
(L w)(n) = a(n) \, w(n+1) + a(n-1) \, w(n-1) + b(n)\, w(n) = \lambda w(n),
\qquad
n \in \mathbb{Z},
\label{B1}
\end{equation}
where $\lambda \in \mathbb{C}$ is the spectral parameter.

Following \cite{D-T} we introduce the two solutions $u(n) = u(n; \lambda)$ and $v(n) = v(n; \lambda)$ of \eqref{B1} which satisfy the initial
conditions
\begin{equation}
u(-1) = 0, \quad u(0) = 1,
\qquad\qquad
v(-1) = -\frac{1}{a(-1)}, \quad v(0) = 0.
\label{B2}
\end{equation}
For $n \geq 0$ the solution $u(n) = u(n; \lambda)$ is a polynomial in $\lambda$ of degree $n$ having the form
\begin{equation}
u(n; \lambda) = \frac{1}{\prod_{j=0}^{n-1} a(j)}
\left\{\lambda^n - \left[\sum_{j=0}^{n-1} b(j)\right] \lambda^{n-1}
+ \left[\sum_{0\leq j < k \leq n-1} b(j) \, b(k) - \sum_{j=0}^{n-2} a(j)^2\right] \lambda^{n-2}
+ \cdots\right\},
\label{DB2}
\end{equation}
while for $n \geq 1$ the solution $v(n) = v(n; \lambda)$ is a polynomial in $\lambda$ of degree $n-1$ having the form
\begin{equation}
v(n; \lambda) = \frac{1}{\prod_{j=0}^{n-1} a(j)}
\left\{\lambda^{n-1} - \left[\sum_{j=1}^{n-1} b(j)\right] \lambda^{n-2}
+ \left[\sum_{1\leq j < k \leq n-1} b(j) \, b(k) - \sum_{j=1}^{n-2} a(j)^2\right] \lambda^{n-3}
+ \cdots\right\},
\label{DB3}
\end{equation}
where we follow the standard convention that empty sums equal $0$, while empty products equal $1$ (e.g., $v(1; \lambda) = a(0)^{-1}$).

Notice also that
\begin{equation}
\left|
  \begin{array}{cc}
    u(n) & \ v(n) \\
    -a(n-1)\,u(n-1) & \ -a(n-1)\,v(n-1) \\
  \end{array}
\right|
= 1
\qquad \text{for all }\;
n \in \mathbb{Z}, \ \lambda \in \mathbb{C}.
\label{B3}
\end{equation}
In particular, $u(n)$ and $v(n)$ are linearly independent solutions of \eqref{B1} for any value of the parameter $\lambda$.

Sometimes it is more convenient, instead of the solutions $u(n)$ and $v(n)$ to work with the (linearly independent) solutions
$\chi(n) = \chi(n; \lambda)$ and $\gamma(n) = \gamma(n; \lambda)$ of \eqref{B1} determined by the initial conditions
\begin{equation}
\chi(0) = 1, \quad \chi(1) = 0
\qquad\text{and}\qquad
\gamma(0) = 0, \quad \gamma(1) = 1.
\label{B2a}
\end{equation}
It follows easily that
\begin{equation}
\chi(n; \lambda) = u(n; \lambda) + [b(0) - \lambda] v(n; \lambda)
\qquad \text{and} \qquad
\gamma(n; \lambda) = a(0) \, v(n; \lambda).
\label{B2b}
\end{equation}
For $n \geq 2$, we have that $\chi(n; \lambda)$ is a polynomial in $\lambda$ of degree $n-2$ and $\gamma(n; \lambda)$ is a polynomial in $\lambda$ of
degree $n-1$. Finally, an easy calculation yields
\begin{equation}
\left|
  \begin{array}{cc}
    \chi(n) & \ \gamma(n) \\
    \chi(n+1) & \ \gamma(n+1) \\
  \end{array}
\right|
= \frac{a(0)}{a(n)}
\qquad \text{for all }\;
n \in \mathbb{Z}, \ \lambda \in \mathbb{C}.
\label{B3a}
\end{equation}

\medskip

\textbf{Remark 1.}  For any fixed $n \geq 0$ the polynomials $v(n; \lambda)$ and $v(n+1; \lambda)$ do not have any common zeros (i.e. they are
relatively prime). The justification of this fact is very simple: Suppose $v(n; \lambda_0) = v(n+1; \lambda_0) = 0$. Then, the fact that
$v(n; \lambda_0)$ satisfies the difference equation \eqref{B1} (for $\lambda = \lambda_0$) implies that $v(n; \lambda_0) = 0$ for all
$n \in \mathbb{Z}$, which is a contradiction since, e.g., $v(1; \lambda_0) = 1/a(0)$. Likewise, the polynomials $u(n; \lambda)$ and $u(n+1; \lambda)$ do not share any common zeros for any fixed $n \geq 0$ and the same is true for $\chi(n; \lambda)$ and $\chi(n+1; \lambda)$ as well as for
$\gamma(n; \lambda)$ and $\gamma(n+1; \lambda)$.
\hfill $\diamondsuit$

\medskip

(The symbol $\diamondsuit$ indicates the end of a remark).

Now let $S$ be the ``$N$-shift" operator
\begin{equation}
(Sf)(n) := f(n + N).
\label{B4}
\end{equation}
Our assumption $a(n+N) = a(n)$ and $b(n+N) = b(n)$ for all $n \in \mathbb{Z}$ implies that the linear operator $S$ maps solutions of \eqref{B1} to
solutions of \eqref{B1} for the same value of $\lambda$ (in other words, $S$ commutes with $L$), and by exploiting this very simple property we can derive the (Floquet) spectral theory of $L$.

For each $\lambda \in \mathbb{C}$ let $\mathcal{W} = \mathcal{W}(\lambda)$
be the two-dimensional vector space of the solutions of \eqref{B1}. By the previous discussion, for each
$\lambda \in \mathbb{C}$ the solutions $u$ and $v$ of \eqref{B1} can be taken
as a basis of $\mathcal{W}(\lambda)$, and the matrix of the operator $S|_{\mathcal{W}}$ with respect to the basis $(u, v)$ is
\begin{equation}
S = S(\lambda) =
\left[
  \begin{array}{cc}
    u(N; \lambda) & \ v(N; \lambda) \\
    -a(-1) \, u(N-1; \lambda) & \ -a(-1) \, v(N-1; \lambda) \\
  \end{array}
\right],
\label{B5}
\end{equation}
where by a little abuse of notation we have also denoted by $S$ the matrix of the operator $S$ restricted to $\mathcal{W}$ (we should not forget that the matrix $S$ and the vector space $\mathcal{W}$ depend on $\lambda$). Thus, $S$ is the Floquet (or monodromy)
matrix associated to equation \eqref{B1}, and formula \eqref{B3} together with the fact that $a(n)$ is $N$-periodic yield immediately that
\begin{equation}
\det S(\lambda) \equiv 1.
\label{B6}
\end{equation}
It follows that the characteristic polynomial of $S(\lambda)$ has the form
\begin{equation}
\det \left(S - r I\right) = r^2 - \Delta(\lambda) \, r + 1,
\quad \text{where}\;
\Delta(\lambda) := \text{tr} S = u(N; \lambda) - a(-1) \, v(N-1; \lambda).
\label{B7}
\end{equation}
Also, it is easy to check that $\Delta(\lambda)$ can be expressed in terms of the solutions $\chi(n; \lambda)$ and $\gamma(n; \lambda)$
(recall \eqref{B2a}) as
\begin{equation}
\Delta(\lambda) =  \chi(N; \lambda) + \gamma(N+1; \lambda).
\label{B7a}
\end{equation}

The quantity $\Delta(\lambda)$ is the (discrete) Hill discriminant of $L$ and from
\eqref{B7}, \eqref{DA3}, \eqref{DB2}, \eqref{DB3}, and \eqref{DA0} it follows that it is a polynomial of $\lambda$  of degree $N$ having the form
\begin{equation}
\Delta(\lambda)
= (-1)^N
\left\{\lambda^N - B_0 N \lambda^{N-1}
+ \left[\sum_{1 \leq j < k \leq N} b(j) \, b(k) - \sum_{j=1}^N a(j)^2\right] \lambda^{N-2}
+ \cdots\right\}.
\label{DB4}
\end{equation}

The eigenvalues $r_1(\lambda)$ and $r_2(\lambda)$ of $S$ are the Floquet multipliers, while their corresponding eigenvectors
$\phi_1(n; \lambda)$ and $\phi_2(n; \lambda)$ are the Floquet solutions of \eqref{B1} so that
\begin{equation}
\phi_j(n + N) = (S \phi_j)(n) = r_j \phi_j(n),
\qquad
j= 1,2.
\label{B8}
\end{equation}

From \eqref{B7} we have
\begin{equation}
r_1(\lambda) \, r_2(\lambda) \equiv 1
\qquad \text{and} \qquad
r_1(\lambda) + r_2(\lambda) = \Delta(\lambda),
\label{B9}
\end{equation}
so that
\begin{equation}
r_1(\lambda), r_2(\lambda) = \frac{\Delta(\lambda) \pm\sqrt{\Delta(\lambda)^2 - 4}}{2}.
\label{B10}
\end{equation}
%
Let us also notice that, by \eqref{B9} $S(\lambda)$ can have a Jordan anomaly only if $r_1(\lambda) = r_2(\lambda) = \pm 1$ (equivalently, only if
$\Delta(\lambda) = \pm 2$) and in the presence of such an anomaly the matrix $S(\lambda)$ is similar to the canonical matrix
\begin{equation*}
\left[
  \begin{array}{cc}
    \pm 1 & 1 \\
    0 & \pm 1 \\
  \end{array}
\right].
\end{equation*}
If this is the case, then there is only one Floquet solution $\phi(n)$ satisfying $\phi(n + N) = \pm \phi(n)$, while there is a second
solution $g(n)$ (sometimes called a generalized Floquet solution), linearly independent to $\phi(n)$, satisfying
\begin{equation}
g(n + N) = \pm g(n) + \phi(n)
\qquad \text{for all }\; n \in \mathbb{Z}.
\label{B13}
\end{equation}
Recall, however, that even if $r_1(\lambda) = r_2(\lambda) = \pm 1$, the Floquet matrix may still be diagonalizable
(and, hence, $S(\lambda) = \pm I$, where $I$ is the $2 \times 2$ identity matrix), in which case we have coexistence of two periodic
(if $S(\lambda) = I$) or antiperiodic (if $S(\lambda) = -I$), linearly independent Floquet solutions.

\medskip

\textbf{Remark 2.} Suppose that a function $\phi(x)$ satisfies
\begin{equation}
\phi(x + N) = r \phi(x)
\qquad \text{for all }\; x \in \mathbb{R},
\label{B8a}
\end{equation}
where $r \ne 0$ is a constant. We write
\begin{equation}
r = e^{\beta N}
\label{B8b}
\end{equation}
and set
\begin{equation}
p(x) := e^{-\beta x} \phi(x).
\label{B8c}
\end{equation}
Then, by using \eqref{B8a} and \eqref{B8b} in \eqref{B8c} we see immediately that $p(x)$ is $N$-periodic and $\phi(x)$ can be written as
\begin{equation}
\phi(x) = e^{\beta x} p(x),
\qquad \text{where }\;
p(x + N) = p(x).
\label{B8d}
\end{equation}
Suppose now that $g(x)$ satisfies
\begin{equation}
g(x + N) = r g(x) + \phi(x)
\qquad \text{for all }\; x \in \mathbb{R},
\label{B13a}
\end{equation}
where $\phi(x)$ satisfies \eqref{B8a}. We set
\begin{equation}
p_1(x) := e^{-\beta x} g(x) - \frac{x}{Nr} \, p(x),
\label{B13b}
\end{equation}
where $p(x)$ is given by \eqref{B8c}. Then, in view of \eqref{B13a}, \eqref{B8b}, and \eqref{B8d} we have
\begin{equation}
p_1(x + N) = e^{-\beta x} e^{-\beta N} [r g(x) + \phi(x)] - \frac{x}{Nr} \, p(x) - \frac{1}{r} \,  p(x) = p_1(x).
\label{B13c}
\end{equation}
Therefore, $g(x)$ can be expressed as
\begin{equation}
g(x) = e^{\beta x} p_1(x) + \frac{x}{Nr} \, e^{\beta x} p(x) = e^{\beta x} p_1(x) + \frac{x}{Nr} \, \phi(x),
\label{B13d}
\end{equation}
where $p_1(x)$ and $p(x)$ are $N$-periodic.

Finally, let us mention that all the above are valid if the functions are defined only for $x \in \mathbb{Z}$, provided, of course,
that $N \in \mathbb{Z}$.
\hfill $\diamondsuit$

\medskip

It is sometimes more convenient to view $r_1(\lambda)$ and $r_2(\lambda)$ as the two branches of a (single-valued) analytic function $r(\lambda)$
defined on the Riemann surface $\Sigma$ of the function $\sqrt{\Delta(\lambda)^2 - 4}$. Then, \eqref{B10} can be
written as
\begin{equation}
r(\lambda) = \frac{\Delta(\lambda) + \sqrt{\Delta(\lambda)^2 - 4}}{2}
\label{B14}
\end{equation}
and $r(\lambda)$ can be called the Floquet multiplier associated to \eqref{B1}.
Let us notice that, since $\Delta(\lambda)^2 - 4$ is a polynomial of even degree, $\Sigma$ has two points at $\infty$. If $\Sigma_{fin}$
denotes the set of
finite points of $\Sigma$ (namely $\Sigma$ minus its two points at $\infty$), then \eqref{B9} implies that $r(\lambda)$ has neither zeros nor poles in $\Sigma_{fin}$. As for the two points at $\infty$ of $\Sigma$, since $\Delta(\lambda)$ has degree $N$ it follows from \eqref{B9} and \eqref{B14}
that at one of these points $r(\lambda)$ has a zero of multiplicity $N$, while at the other it has a pole of order $N$. Also, it maybe worth mentioning (i) that
\eqref{B14} is equivalent to
\begin{equation}
\Delta(\lambda) = r(\lambda) + \frac{1}{r(\lambda)}
\label{B14a}
\end{equation}
and (ii) that
\begin{equation}
\frac{r'(\lambda)}{r(\lambda)} = \frac{\Delta'(\lambda)}{\sqrt{\Delta(\lambda)^2 - 4}}.
\label{B14b}
\end{equation}
The Floquet solutions too can be viewed as the two branches of a meromorphic function defined on $\Sigma$. First we normalize them so that
$\phi_1(0; \lambda) = \phi_2(0; \lambda) = 1$. It, then, follows that
$\phi_1(n; \lambda)$ and $\phi_2(n; \lambda)$ are the branches of the function
\begin{equation}
\phi(n; \lambda) = u(n; \lambda) - \frac{u(N; \lambda) - r(\lambda)}{v(N; \lambda)} \, v(n; \lambda).
\label{B14c}
\end{equation}
As we have already mention, deg$_{\lambda}v(N; \lambda) = N-1$ (where deg$_{\lambda}v(N; \lambda)$ denotes the degree of $v(N; \lambda)$ viewed as
a polynomial of $\lambda$). Hence, $\phi(n; \lambda)$, as a function of $\lambda$, can have at most $N-1$ poles
in $\Sigma_{fin}$ counting multiplicities (in the non-self-adjoint case the zeros of $v(N; \lambda)$  are not necessarily simple ---
an example where $v(4; \lambda)$ has a triple zero is mentioned in Subsection 2.3).

Having $\Delta(\lambda)$ and $r(\lambda)$, the $l^2(\mathbb{Z})$-spectrum $\sigma(L)$ of $L$ can be characterized as
\begin{equation}
\sigma(L) = \{\lambda \in \mathbb{C} \,:\, \Delta(\lambda) \in [-2, 2\,]\}
\quad \Leftrightarrow \quad
\sigma(L) = \{\lambda \in \mathbb{C} \,:\, |r(\lambda)| = 1\},
\label{B15}
\end{equation}
which implies that $\sigma(L)$ is a finite union of bounded analytic arcs lying in the complex plane
(notice that by the first equality in \eqref{B9} we have that $|r_1(\lambda)| = 1$ if and only if $|r_2(\lambda)| = 1$).

The adjoint operator $L^{\ast}$ of $L$ is given by the formula
\begin{equation}
(L^{\ast} w)(x) = \overline{a(n)} \, w(n+1) + \overline{a(n-1)} \, w(n-1) + \overline{b(n)}\, w(n),
\qquad
n \in \mathbb{Z},
\label{B16}
\end{equation}
where the bar denotes complex conjugation. Hence $L$ is self-adjoint if and only if $a(n)$ and $b(n)$ are real-valued. In general we have
\begin{equation}
\sigma(L^{\ast}) = \overline{\sigma(L)},
\label{B17}
\end{equation}
namely $\lambda \in \sigma(L^{\ast})$ if and only if $\overline{\lambda} \in \sigma(L)$. If in particular $\sigma(L) \subset \mathbb{R}$, then $\sigma(L^{\ast}) = \sigma(L)$.

\subsection{Floquet spectrum; periodic and antiperiodic eigenvalues}

In view of \eqref{B9} and \eqref{B10}, formula \eqref{B15} yields some additional characterizations of the $l^2(\mathbb{Z})$-spectrum of $L$,
namely
\begin{equation}
\sigma(L) = \{\lambda \in \mathbb{C} \,:\, \Delta(\lambda) = 2 \cos(\kappa), \ \; 0 \leq \kappa \leq \pi\}
\label{K1a}
\end{equation}
and
\begin{equation}
\sigma(L) = \{\lambda \in \mathbb{C} \,:\, r_1(\lambda) = e^{i \kappa}, \ \; 0 \leq \kappa \leq \pi\}.
\label{K1b}
\end{equation}
Thus, if for a given $\kappa \in [0, \pi]$, we introduce the \emph{Floquet spectrum}
\begin{equation}
\sigma_{\kappa}(L) := \{\lambda \in \mathbb{C} \,:\, \Delta(\lambda) = 2 \cos(\kappa)\},
\label{K2}
\end{equation}
then $\sigma(L)$ can be written as the disjoint union
\begin{equation}
\sigma(L) = \bigcup_{0 \leq \kappa \leq \pi} \sigma_{\kappa}(L).
\label{K3}
\end{equation}
Clearly, the Floquet spectrum $\sigma_{\kappa}(L)$ is the set of zeros of the $N$-th degree polynomial
\begin{equation}
F_{\kappa}(\lambda) := \Delta(\lambda) - 2 \cos(\kappa).
\label{K2a}
\end{equation}
Observe that $F_{\kappa}'(\lambda) = \Delta'(\lambda)$ is independent of $\kappa$ and has degree $N-1$. Thus, if $\lambda$ is a multiple zero of
$F_{\kappa}(\lambda)$, then $\lambda$ must be a zero of $\Delta'(\lambda)$, and there are at most $N-1$ such zeros (which, of course, are independent of $\kappa$). For each such value of $\lambda$ there is at most one $\kappa \in [0, \pi]$ for which $F_{\kappa}(\lambda) = 0$ (since $\cos(\kappa)$
is strictly decreasing on $[0, \pi]$). It follows that there are at most $N-1$ values of $\kappa \in [0, \pi]$ for which $F_{\kappa}(\lambda)$ has
multiple zeros and, therefore, if $\kappa$ is not equal to any of those exceptional values, the Floquet spectrum $\sigma_{\kappa}(L)$ consists of $N$ distinct $\kappa$\emph{-Floquet eigenvalues}. Recall, e.g., that in the self-adjoint case, if $\kappa \ne 0, \pi$, then $F_{\kappa}(\lambda)$ has $N$
distinct zeros.

Let us first consider the case $\kappa \in (0, \pi)$, namely $\kappa \ne 0$ and $\kappa \ne \pi$. Under this assumption for $\kappa$,
if $\lambda \in \sigma_{\kappa}(L)$, we have $r_1(\lambda) = e^{i \kappa} \ne \pm 1$ and, therefore, $\lambda$ is not a branch point of $r(\lambda)$,
hence there are two linearly independent Floquet solutions $\phi_1(n) = \phi_1(n; \lambda)$ and $\phi_1(n) = \phi_2(n; \lambda)$ corresponding to any particular $\lambda \in \sigma_{\kappa}(L)$, satisfying
\begin{equation}
\phi_1(n + N) = e^{i \kappa} \phi_1(n),
\qquad
n \in \mathbb{Z}
\label{K4a}
\end{equation}
and
\begin{equation}
\phi_2(n + N) = e^{-i \kappa} \phi_2(n),
\qquad
n \in \mathbb{Z}.
\label{K4b}
\end{equation}
In view of \eqref{B8d} of Remark 2, formula \eqref{K4a} implies that
\begin{equation}
\phi_1(n) = e^{i \kappa n/N} p(n),
\qquad \text{where }\;
p(n + N) = p(n).
\label{NF1}
\end{equation}
By using \eqref{NF1} in \eqref{B1} we can see that $p(n)$ satisfies the boundary value problem
\begin{align}
a(n) e^{i \kappa /N} p(n+1) + a(n-1) e^{-i \kappa /N} p(n-1) &+ b(n)\, p(n) = \lambda p(n),
\nonumber
\\
&\ \label{NF2}
\\
p(0) = p(N),
\qquad\qquad
p(1) &= p(N+1).
\nonumber
\end{align}
Problem \eqref{NF2} can be written in the matrix form
\begin{equation}
M_{\kappa} \vec{p} = \lambda \vec{p}
\label{NF3}
\end{equation}
where $\vec{p}$ is the column vector $\vec{p} := [p(0), \ldots, p(N-1)]^{\top}$ and $M_{\kappa}$ is the $N \times N$ matrix (for $N \geq 3$)
\begin{equation}
M_{\kappa} :=
\left[
  \begin{array}{cccccc}
    b(0) & a(0)e^{i \kappa /N} & 0 & \cdots & 0 & a(N-1)e^{-i \kappa /N} \\
    a(0)e^{-i \kappa /N} & b(1) & a(1)e^{i \kappa /N} & \cdots & 0 & 0 \\
		0 & a(1)e^{-i \kappa /N} & b(2) & \cdots & 0 & 0 \\
		\vdots & \vdots & \vdots & \ddots & \vdots & \vdots \\
		0 & 0 & 0 & \cdots & b(N-2) & a(N-2)e^{i \kappa /N} \\
		a(N-1)e^{i \kappa /N} & 0 & 0 & \cdots & a(N-2)e^{-i \kappa /N} & b(N-1) \\
  \end{array}
\right].
\label{Nf4}
\end{equation}
If $\kappa$ is such that the polynomial $F_{\kappa}(\lambda)$ of \eqref{K2a} has simple zeros, then $F_{\kappa}(\lambda)$ must be the characteristic
polynomial of $M_{\kappa}$. Then, by continuity we have that
\begin{equation}
\det\left(M_{\kappa} - \lambda I\right) = \Delta(\lambda) - 2 \cos(\kappa)
\qquad \text{for all }\;
\kappa \in [0, \pi]
\label{NF5}
\end{equation}
(and, consequently, by analytic continuation the above equation must hold for all $\kappa \in \mathbb{C}$). In particular, the spectrum of
$M_{\kappa}$ is $\sigma_{\kappa}(L)$. Furthermore, the (pure) eigenvectors of $M_{\kappa}$ correspond precisely to the Floquet solutions satisfying
\eqref{K4a}.

If $\lambda \in \mathbb{C}$ is such that one of the Floquet multipliers is equal to $1$, then $\lambda$ is a \emph{periodic eigenvalue} of $L$.
In this case \eqref{B9} and \eqref{B14} imply that $r_1(\lambda) = r_2(\lambda) = r(\lambda) = 1$ or, equivalently, $\Delta(\lambda) = 2$.
Likewise, if $\lambda \in \mathbb{C}$ is such that one of the Floquet multipliers is equal to $-1$, then $\lambda$ is an
\emph{antiperiodic eigenvalue} of $L$.
In this case \eqref{B9} and \eqref{B14} imply that $r_1(\lambda) = r_2(\lambda) = r(\lambda) = -1$ or, equivalently, $\Delta(\lambda) = -2$.
Furthermore, by \eqref{B15} we have that periodic and antiperiodic eigenvalues are always in the spectrum $\sigma(L)$ of $L$.

If $\lambda \in \mathbb{C}$ is neither a periodic nor an antiperiodic eigenvalue of $L$, then from \eqref{B9} we get that
$r_1(\lambda) \ne r_2(\lambda)$, which in turn implies that $\lambda$ is not a branch point of $r(\lambda)$ (and, also, as we have already seen,
that the Floquet matrix $S(\lambda)$ is diagonalizable). Thus, branch points of $r(\lambda)$ as well as Jordan anomalies of $S(\lambda)$ can occur
only at periodic or antiperiodic eigenvalues of $L$, and for these reasons these eigenvalues are quite special. Let us recall that in the self-adjoint case $\lambda$ is a branch point of $r(\lambda)$ if and only if $S(\lambda)$ has a Jordan anomaly (and such a $\lambda$ must necessarily be real).
However, this is not always true in the non-self-adjoint case. We will encounter this phenomenon in Subsection 4.2.

If $w(n)$ is a periodic or antiperiodic (Floquet) solution of \eqref{B1} (for a fixed $\lambda$), then, obviously, $w(n + 2N) = u(n)$, i.e. $w(n)$ is
$2N$-periodic. Conversely, if for a fixed $\lambda$ the equation \eqref{B1} posesses a $2N$-periodic solution $w(n)$, then, in view of \eqref{B4}, we get $(S^{\,2} w)(n) = w(n + 2N) = w(n)$. Thus, the operator $(S|_{\mathcal{W}})^2$, acting of the two-dimensional space
$\mathcal{W} = \mathcal{W}(\lambda)$, has an eigenvalue equal to $1$, and, consequently, both eigenvalues $r_1(\lambda)^2$ and $r_2(\lambda)^2$ of
$(S|_{\mathcal{W}})^2$ are equal to $1$. Hence $r_1(\lambda) = r_2(\lambda) = \pm1$ and
$w(n)$ is a periodic or an antiperiodic (Floquet) solution of \eqref{B1}.

Now, for an integer $m \geq 1$ let us consider the space
\begin{equation}
\mathcal{P}_m := \{f(n) : f(n + m) = f(n) \ \ \ \text{for all }\, n \in \mathbb{Z}\},
\label{G1}
\end{equation}
namely the set of $m$-periodic sequences over the complex numbers.
Obviously, $\mathcal{P}_m$ is a vector space of (complex) dimension $m$. In the case where $m$ is a multiple of $N$, the operator $L$, having $N$-periodic coefficients $a(n)$ and $b(n)$, maps $\mathcal{P}_m$ into $\mathcal{P}_m$. In particular, for $m = 2N$ the operator $L$ maps
$\mathcal{P}_{2N}$ into $\mathcal{P}_{2N}$ and, due to the previous discussion this is the most interesting case.
As a basis of $\mathcal{P}_{2N}$ we can choose the sequences
\begin{equation}
e_j(n) := \delta_{jn},
\qquad
n \in \mathbb{Z},\ \; j = 1, \dots, 2N,
\label{G2}
\end{equation}
where $\delta_{jn}$ is the Kronecker delta. Then, the $2N \times 2N$ matrix of $L|_{\mathcal{P}_{2N}}$ with respect to that basis is
\begin{equation}
L_{2N} :=
\left[
  \begin{array}{ccccccc}
    b(1) & a(1) & 0 & \cdots & 0 & 0 & a(2N) \\
    a(1) & b(2) & a(2) & \cdots & 0 & 0 & 0 \\
		0 & a(2) & b(3) & \cdots & 0 & 0 & 0 \\
		\vdots & \vdots & \vdots & \ddots & \vdots & \vdots & \vdots \\
		0 & 0 & 0 & \cdots & b(2N-2) & a(2N-2) & 0 \\
		0 & 0 & 0 & \cdots & a(2N-2) & b(2N-1) & a(2N-1) \\
		a(2N) & 0 & 0 & \cdots & 0 & a(2N-1) & b(2N) \\
  \end{array}
\right],
\label{G3}
\end{equation}
where $a(n + N) = a(n)$ and $b(n + N) = b(n)$.
Notice that the matrix $L_{2N}$ is symmetric, but not Hermitian, unless, of course $a(n)$ and $b(n)$ are real-valued, in which case $L_{2N}$ is real symmetric (hence Hermitian) and its associated operator $L|_{\mathcal{P}_{2N}}$ is self-adjoint.

It follows that the eigenvectors of $L|_{\mathcal{P}_{2N}}$ (in $\mathcal{P}_{2N}$), being $2N$-periodic solutions of \eqref{B1}, are precisely the $N$-periodic and $N$-antiperiodic (linearly independent) solutions of \eqref{B1}. Also, the spectrum of the operator $L|_{\mathcal{P}_{2N}}$,
i.e. the set of eigenvalues of the matrix $L_{2N}$, coincides with the set of periodic and antiperiodic eigenvalues of $L$, that is the zeros of the
polynomials $\Delta(\lambda) - 2$ and $\Delta(\lambda) + 2$. From \eqref{DB4} we know that $\Delta(\lambda)$ has degree $N$. Hence the polynomial
$\Delta(\lambda) - 2$ has at most $N$ distinct zeros and the same is true for $\Delta(\lambda) + 2$. Obviously, these two polynomials cannot have
common zeros. On the other hand $\Delta(\lambda) - 2$ and $\Delta(\lambda) + 2$ have the same derivative, namely $\Delta'(\lambda)$, from which it
follows that $\Delta(\lambda)^2 - 4$ has at least $N + 1$ distinct zeros (the derivative of $\Delta(\lambda)^2 - 4$ is
$2 \Delta(\lambda) \Delta'(\lambda)$ and $\Delta(\lambda)$ does not have common zeros with $\Delta(\lambda)^2 - 4$).

From \eqref{DB4} we get that $\Delta(\lambda)^2 - 4$ is a monic polynomial of degree $2N$, i.e. its leading term is $\lambda^{2N}$. Also, for
generic $a(n)$ and $b(n)$ we have that $\Delta(\lambda)^2 - 4$ has $2N$ simple zeros. Therefore, the characteristic polynomial of $L_{2N}$ is
\begin{equation}
\det\left(L_{2N} - \lambda I\right) = \Delta(\lambda)^2 - 4.
\label{NF6}
\end{equation}
Another way to justify \eqref{NF6} is by using an argument similar to one used for establishing \eqref{NF5}.

\subsection{Certain classes of isospectral operators}
Here we will present some classes of operators sharing the same spectral properties.

Let $a^l(n) := a(n+l)$ and $b^{\,l}(n) := b(n+l)$, where $l \in \mathbb{Z}$, and consider the operator
\begin{equation}
(L^l w)(n) := a^l(n) \, w(n+1) + a^l(n-1) \, w(n-1) + b^{\,l}(n)\, w(n),
\qquad
n \in \mathbb{Z}.
\label{B15a}
\end{equation}
Then, by writing equation \eqref{B8} as
\begin{equation}
\phi_j(n  + l + N; \lambda) = r_j(\lambda) \phi_j(n + l; \lambda),
\qquad
j= 1,2.
\label{B15b}
\end{equation}
we can see that $L$ and $L^l$ have the same multiplier $r(\lambda)$ and, consequently, the same discriminant.
Hence, $\sigma(L^l) = \sigma(L)$. Also, if $a^{\sharp}(n) := a(-n)$, $b^{\,\sharp}(n) := b(-n)$, and $L^{\sharp}$ is the operator associated to $a^{\sharp}$ and $b^{\,\sharp}$,
then $\phi(-n; \lambda)$ is a Floquet solution of $L^{\sharp} w = \lambda w$ if and only if $\phi(n; \lambda)$ is a Floquet solution of
$L w = \lambda w$. It follows that $L$ and $L^{\sharp}$ have the same multiplier, the same discriminant, and the same spectrum.

This is to be compared with the continuous case where the potentials
$V(x)$, $V^{\xi}(x) := V(x + \xi)$, and $V^{\sharp}(x) := V(-x)$ have the same multiplier (hence the same discriminant and the same
$L^2(\mathbb{R})$-spectrum) for all $\xi \in \mathbb{R}$) even for nonreal values of $\xi$, as long as
$V(x)$ is analytic and $V(x + \xi)$ makes sense.

We continue with another case of isospectrality. Formula \eqref{DB4} hints that the discriminant $\Delta(\lambda)$ may stay unchanged if for some $n$
we replace $a(n)$ by $-a(n)$. Actually this guess is essentially true as it follows from a simple fact mentioned in \cite{T} for the self-adjoint
case, which extends automatically in the case of complex $a(n)$ and $b(n)$.

\medskip

\textbf{Proposition 1.} Suppose $\tau(n) = 1$ or $-1$ and $\tau(n+N) = \tau(n)$ for all $n \in \mathbb{Z}$. Let $a(n)$ and
$b(n)$ be the coefficients of the ($N$-periodic) Jacobi operator $L$ of \eqref{A3} and consider the operator $\hat{L}$ whose coefficients are
$\hat{a}(n) := \tau(n) a(n)$ and $\hat{b}(n) := b(n)$. If $\Delta(\lambda)$ and $\hat{\Delta}(\lambda)$ are the discriminants of $L$ and
$\hat{L}$ respectively, then
\begin{equation}
\hat{\Delta}(\lambda) = (-1)^{\nu(\tau)} \Delta(\lambda),
\qquad \text{where }\
\nu(\tau) := \#\{n \, : \, \tau(n) = -1, \ 0 \leq n \leq N-1\}
\label{I1}
\end{equation}
(as usual, $\# \mathcal{S}$ denotes the cardinality of the set $\mathcal{S}$). In particular, $\sigma(\hat{L}) = \sigma(L)$.

\smallskip

\textit{Proof}. Suppose $w(n)$ satisfies $(L w)(n) = \lambda w(n)$. Then it is easy to see that
\begin{equation}
\hat{w}(n) :=
\left\{
\begin{array}{ccc}
    w(n) \prod_{j=0}^{n-1} \tau(j), & \qquad & n > 0, \\
    w(0), & \qquad & n = 0, \\
		w(n) \prod_{j=n}^{-1} \tau(j), & \qquad & n < 0, \\
  \end{array}
\right.
\label{I2}
\end{equation}
satisfies
$(\hat{L} \hat{w})(n) = \lambda \hat{w}(n)$. In particular, $w(n)$ is a Floquet solution associated to $L$, with multiplier $r$, if and only if
$\hat{w}(n)$ is a Floquet solution associated to $\hat{L}$ with multiplier $\hat{r} := (-1)^{\nu(\tau)} r$, and the proof is completed by recalling that $\Delta(\lambda) = r(\lambda) + r(\lambda)^{-1}$ and $\hat{\Delta}(\lambda) = \hat{r}(\lambda) + \hat{r}(\lambda)^{-1}$.
\hfill $\blacksquare$

\medskip

(The symbol $\blacksquare$ indicates the end of a proof).

Motivated by Proposition 1 we introduce the following equivalence relation between Jacobi operators.

\medskip

\textbf{Definition 1.} Two Jacobi operators $L$ and $\hat{L}$ (of complex coefficients) are called equivalent, symbolically $L \sim \hat{L}$, if
their associated coefficients $a(n)$, $b(n)$, $\hat{a}(n)$, and $\hat{b}(n)$ are related as $\hat{a}(n) = \tau(n) a(n)$ and $\hat{b}(n) = b(n)$
for all $n \in \mathbb{Z}$, where $\tau(n) = 1$ or $-1$.

\medskip

In other words, $L \sim \hat{L}$, if $a(n)^2 = \hat{a}(n)^2$ and $\hat{b}(n) = b(n)$ for all $n \in \mathbb{Z}$.

\medskip

\textbf{Remark 3.} Clearly, in the $N$-periodic case, if the coefficient $a(n)$ of $L$ satisfies the normalization \eqref{DA0}, then the coefficient
$\hat{a}(n) = \tau(n) a(n)$ of $\hat{L}$ satisfies \eqref{DA0} if and only if ${\nu(\tau)}$ of \eqref{I1} satisfies $(-1)^{\nu(\tau)} = 1$ (i.e. ${\nu(\tau)}$ is even). Thus, if the coefficients of both $L$ and $\hat{L}$ satisfy the normalization \eqref{DA0}
and $L \sim \hat{L}$, then $\hat{\Delta}(\lambda) = \Delta(\lambda)$.
\hfill $\diamondsuit$

\subsection{The Dirichlet spectrum}
Let us look at the Dirichlet-type boundary value problem ($N \geq 2$)
\begin{align}
(L \psi)(n) &= a(n) \, \psi(n+1) + a(n-1) \, \psi(n-1) + b(n)\, \psi(n) = \mu \, \psi(n)
\label{D1a}
\\
\psi(0) &= \psi(N) = 0
\label{D1b}
\end{align}
(notice that $\psi(n)$ can be extended so that it satisfies \eqref{D1a} for all $n \in \mathbb{Z}$).
Clearly, the eigenvalues of the problem \eqref{D1a}--\eqref{D1b} are the zeros of the polynomial $v(N; \lambda)$. As we have seen
$\text{deg}_{\lambda}v(N; \lambda) = N-1$, hence there are $N-1$ Dirichlet eigenvalues $\mu_1, \dots, \mu_{N-1}$, counting multiplicities.
Hence, in view of \eqref{DB3}, \eqref{DA0}, and \eqref{B2b}
\begin{equation}
v(N; \lambda) = (-1)^N \prod_{j=1}^{N-1} \left(\lambda - \mu_j\right)
\qquad \text{and} \qquad
\gamma(N; \lambda) = (-1)^N a(0)\prod_{j=1}^{N-1} \left(\lambda - \mu_j\right)
\label{D2}
\end{equation}

In the case where $a(n)$ and $b(n)$ are real-valued, the problem \eqref{D1a}--\eqref{D1b} is self-adjoint and hence the eigenfunctions form a basis
of the underlying vector space, which is clearly $(N-1)$-dimensional. Since for each $\mu_j$ we cannot have more than one Dirichlet eigenfunction
(up to linear independence), it follows that in the real case the zeros of $v(N; \lambda)$ are real and simple (and between any two bands of the
spectrum there is exactly one Dirichlet eigenvalue). However, this is not always true in the case of nonreal $a(n)$, $b(n)$. For example, if $N=4$
and $a(n) \equiv -1$, then
\begin{equation*}
v(4; \lambda) = -[\lambda - b(1)] [\lambda - b(2)] [\lambda - b(3)]) + \lambda - b(1) + \lambda - b(3),
\end{equation*}
and the choice $b(1) = -b(3) = \sqrt{2} \, i$ and $b(2) = 0$ yields $v(4; \lambda) = -\lambda^3$, hence  $\mu_1 = \mu_2 = \mu_3 = 0$. Fixing
$b(4) = 0$ gives a specific ``pathological" example, namely
\begin{equation*}
a(n) \equiv -1,
\qquad
b(n) = \frac{i^n - (-i)^n}{\sqrt{2}}.
\end{equation*}


Next, let us observe that by invoking \eqref{DB3} we get immediately the ``trace formula"
\begin{equation}
\mu_1 + \cdots + \mu_{N-1} = b(1) + \cdots + b(N-1),
\label{D4a}
\end{equation}
which, in view of \eqref{DA3} can be also written as
\begin{equation}
\mu_1 + \cdots + \mu_{N-1} = B_0 N - b(0).
\label{D4b}
\end{equation}
Also, since by \eqref{DA3} and \eqref{DB4} we have
\begin{equation}
\Delta(\lambda)^2 - 4
= \left\{\lambda^{2N} - 2 \left[b(1) + \cdots + b(N)\right] \lambda^{2N-1}
+ \cdots\right\},
\label{D4c}
\end{equation}
it follows that (recall \eqref{DA3})
\begin{equation}
\sum_{j=0}^{2N-1}\lambda_j = 2 \left[b(1) + \cdots + b(N)\right] = 2 B_0 N,
\label{D4e}
\end{equation}
where $\lambda_j$, $j = 0, 1, \dots, 2N-1$, are the zeros of $\Delta(\lambda)^2 - 4$ (counting multiplicities), namely the periodic and
antiperiodic eigenvalues. Furthermore, \eqref{D4b} yields
\begin{equation}
\sum_{j=0}^{2N-1}\lambda_j - 2\sum_{j=1}^{N-1} \mu_j = 2b(0).
\label{D4d}
\end{equation}

Finally, let us mention that the Dirichlet eigenfunction $\psi(n)$, extended to $\mathbb{Z}$, is always a
Floquet solution.

\subsection{The unperturbed case}
If $a(n) \equiv -1$ and $b(n) \equiv 0$ (viewed as $N$-periodic functions), then the operator $L$ reduces to the unperturbed operator
\begin{equation}
\left(\tilde{L} w\right)(n) := -w(n+1) - w(n-1),
\qquad
n \in \mathbb{Z},
\label{DC0}
\end{equation}
and equation \eqref{B1} becomes
\begin{equation}
\left(\tilde{L} w\right)(n) = -w(n+1) - w(n-1) = \lambda w(n),
\qquad
n \in \mathbb{Z}.
\label{DC0a}
\end{equation}
From now on a tilded quantity will be always associated with the unperturbed case.

It is convenient to introduce a new spectral parameter $z$ related to $\lambda$ as
\begin{equation}
z + z^{-1} := -\lambda.
\label{DC0b}
\end{equation}
Then, the solutions $\chi$ and $\gamma$ (recall \eqref{B2a}) in the unperturbed case become respectively
\begin{equation}
\tilde{\chi}(n; \lambda) = \frac{z^{1-n} - z^{n-1}}{z - z^{-1}},
\qquad \text{and} \qquad
\tilde{\gamma}(n; \lambda) = \frac{z^n - z^{-n}}{z - z^{-1}}.
\label{DC0c}
\end{equation}
In particular, for $\lambda = -2$ (equivalently $z = 1$) we have
\begin{equation}
\tilde{\chi}(n; -2) = 1-n
\qquad \text{and} \qquad
\tilde{\gamma}(n; -2) = n,
\label{DC0cc}
\end{equation}
while for $\lambda = 2$ (equivalently $z = -1$) we have
\begin{equation}
\tilde{\chi}(n; 2) = (-1)^{n-1} n
\qquad \text{and} \qquad
\tilde{\gamma}(n; 2) = (-1)^{n-1} (n-1).
\label{DC0ccc}
\end{equation}
By straightforward induction we can also see that the solution $\tilde{\gamma}(n; \lambda)$, $n \geq 3$, expanded in
descending powers of $\lambda$, has the form
\begin{equation}
\tilde{\gamma}(n; \lambda) = (-1)^{n-1} \lambda^{n-1} + (-1)^n (n-2) \lambda^{n-3} + \cdots
\label{DC0c7}
\end{equation}
Now, using \eqref{DC0c} in \eqref{B7a} we get that the discriminant of the unperturbed operator is
\begin{equation}
\tilde{\Delta}_N(\lambda) := z^N + z^{-N}
= \left(\frac{-\lambda +\sqrt{\lambda^2 - 4}}{2}\right)^N + \left(\frac{-\lambda -\sqrt{\lambda^2 - 4}}{2}\right)^N.
\label{DC0c4}
\end{equation}
Also, from \eqref{DC0c7} and the fact that for $n \geq 2$ we have (in the unperturbed case) that
$\tilde{u}(n; \lambda) = -\tilde{v}(n-1; \lambda)$, we obtain the expansion
\begin{equation}
\tilde{\Delta}_N(\lambda) = (-1)^N \lambda^N - (-1)^N N \lambda^{N-2} + \cdots,
\qquad \text{for }\; N \geq 2
\label{DC0c5}
\end{equation}
(this also follows immediately from \eqref{DB4}). Thus, in particular, the coefficient of $\lambda^{N-1}$ in $\tilde{\Delta}_N(\lambda)$ is $0$,
which can be also seen from \eqref{DB4}.

From the above formulas it follows easily that in the unperturbed case the Floquet multiplier becomes
\begin{equation}
\tilde{r}(\lambda) = z^N = \left(\frac{-\lambda + \sqrt{\lambda^2 - 4}}{2}\right)^N,
\label{DC0c6}
\end{equation}
while the spectrum is
\begin{equation}
\sigma\left(\tilde{L}\right) = [-2, 2].
\label{DC0d}
\end{equation}
Furthermore, if we set
\begin{equation}
z_k := e^{i\pi k/N},
\qquad
k = 0, 1, \ldots, 2N-1,
\label{DC0e}
\end{equation}
and
\begin{equation}
\tilde{\lambda}_k := -\left(z_k + z_k^{-1}\right) = -2 \cos\left(\frac{\pi k}{N}\right),
\qquad
k = 0, 1, \ldots, 2N-1,
\label{DC0f}
\end{equation}
then $\tilde{\lambda}_0 = -2$, $\tilde{\lambda}_N = 2$ and $\tilde{\lambda}_k = \tilde{\lambda}_{2N - k}$ for $k = 1, \ldots, N-1$.
In addition, $\tilde{\lambda}_0$ and $\tilde{\lambda}_N$
are simple zeros of $\tilde{\Delta}_N(\lambda)^2 - 4$, while $\tilde{\lambda}_k$, $k = 1, \ldots, N-1$, are double zeros of
$\tilde{\Delta}_N(\lambda)^2 - 4$.
It follows that $\tilde{\lambda}_0 = -2$ is a periodic eigenvalue of $\tilde{L}$ of geometric multiplicity $1$,
the corresponding eigenfunction being $\tilde{\phi}(n; -2) \equiv 1$, while $\tilde{\lambda}_N = 2$ is a periodic
(antiperiodic) eigenvalue of $\tilde{L}$ of geometric multiplicity $1$, if $N$ is even (odd),
the corresponding eigenfunction being $\tilde{\phi}(n; 2) = (-1)^n$. Finally, each
$\tilde{\lambda}_k = 2\cos(\pi k / N)$, $k = 1, \ldots, N-1$, is a periodic (antiperiodic)
eigenvalue of $\tilde{L}$ of geometric multiplicity $2$, if $k$ is even (odd), while the associated eigenfunctions are
$\tilde{\phi}_1(n; \tilde{\lambda}_k) = e^{i\pi k n/N}$ and $\tilde{\phi}(n; \tilde{\lambda}_k) = e^{-i\pi k n/N}$
(i.e. we have coexistence of two linearly independent periodic or antiperiodic solutions).

\subsection{The essentially unperturbed operators}

\textbf{Definition 2.} We say that $L$ of \eqref{A3} is an \textit{essentially unperturbed operator} if $L \sim \tilde{L}$, i.e. if $a(n)^2 \equiv 1$
and $b(n) \equiv 0$.

\medskip

From the above definition it follows that, for a given period $N$ there are $2^N$ essentially unperturbed operators, one of them being $\tilde{L}$.
Obviously the essentially unperturbed operators have real coefficients and hence they are self-adjoint.
Notice also that $L$ is essentially unperturbed if and only if $-L$ is essentially unperturbed.
If $N$ is odd and $L$ is essentially unperturbed, then either $L$ or $-L$ satisfies the normalization \eqref{DA0}.

\medskip

\textbf{Remark 4.} There are many results which can be proved by first checking that they are valid for the essentially unperturbed case and then view the general case as a continuous deformation of the  unperturbed case. For instance, let us show that for any real $a(n)$ and $b(n)$
(with $a(n) \ne 0$ for all $n$) the zeros of the polynomials
$v(N; \lambda)$ and $v(N+1; \lambda)$ interlace. In the case of $a(n) \equiv \pm1$ and $b(n) \equiv 0$ the statement follows easily from \eqref{DC0b}
and \eqref{DC0c}. Now, given any $a(n) \ne 0$ and $b(n)$ consider the family of quantities $a(n; t) \ne 0$ and $b(n; t)$, $t \in [0, 1]$
continuous in $t$, such that $a(n; 0) = \text{sgn} [a(n)]$, $a(n; 1) = a(n)$, $b(n; 0) = 0$, and $b(n; 1) = b(n)$
(e.g., $b(n; t) = tb(n)$). For each $t$ the zeros of
$v(N; \lambda; t)$ and $v(N+1; \lambda; t)$ (where $v(n; \lambda; t)$ denotes the solution of \eqref{A3} when the coefficients of $L$ are
$a(n; t)$ and $b(n; t)$), such that $v(0; \lambda; t) = 0$ and $v(1; \lambda; t) = 1$) are real, being the Dirichlet eigenvalues of a self-adjoint
operator. Furthermore,
as $t$ moves continuously from $0$ to $1$ no zero of $v(N; \lambda; t)$ can ``cross" a zero of $v(N+1; \lambda; t)$ due to Remark 1. Hence,
the relative position of the zeros of $v(N; \lambda; t)$ and $v(N+1; \lambda; t)$ is independent of $t$. Since for $t=0$ their zeros interlace,
it follows that they also interlace for $t=1$.
\hfill $\diamondsuit$

\section{Inverse spectral considerations for the Schr\"{o}dinger case}

In this section we consider certain inverse spectral problems for the discrete periodic Schr\"{o}dinger operator with a (complex) potential $b(n)$,
namely the operator
\begin{equation}
(L_{\text{Schr}} \, w)(n) := -w(n+1) - w(n-1) + b(n)\, w(n),
\qquad
n \in \mathbb{Z}.
\label{L0}
\end{equation}

\medskip

\textbf{Proposition 2.} For the case where $L = L_{\text{Schr}}$ the zeros of $v(N; \lambda)$ (counting multiplicities) together with the zeros of
$v(N+1; \lambda)$ determine the sequence $b(n)$.

\smallskip

\textit{Proof}. As we have seen, the polynomial $v(N; \lambda)$ has $N-1$ zeros counting multiplicities. Hence, from its zeros we also know $N$.
Then, in view of \eqref{DB3} we know $\sum_{j=1}^{N-1} b(j)$. Likewise, from the zeros of $v(N+1; \lambda)$ we can recover $\sum_{j=1}^N b(n)$.
Hence, from the given data we can get $b(N)$. Having $b(N)$ we can use the difference equation \eqref{B1}, satisfied by $v(n; \lambda)$, in order to
recover $v(N-1; \lambda)$. Having now $v(N; \lambda)$ and $v(N-1; \lambda)$ we can recover $b(N-1)$ and $v(N-2; \lambda)$. We continue in the same
manner until we recover $b(j)$ for all $j = 1, \dots, N$.
\hfill $\blacksquare$

\medskip

An essentially equivalent version of Proposition 2 has appeared in \cite{P-D}. The proposition can be viewed as a special case of a discrete
counterpart of a general result of Levitan and Gasymov \cite{L-G}, in the continuous case, which says that a potential can be recovered from two
spectra.

Let us, now, discuss a variance of the case of the above proposition. Suppose we know the zeros (counting multiplicities) of the polynomials
$v(N; \lambda)$ and $u(N+1; \lambda)$. Then, as we have seen the polynomials $v(N; \lambda)$ and $u(N+1; \lambda)$ are known, while from \eqref{B3} we get
\begin{equation}
u(N; \lambda) \, v(N+1; \lambda) = 1 + v(N; \lambda) \, u(N+1; \lambda),
\label{H1}
\end{equation}
hence the polynomial $F(\lambda) := u(N; \lambda) \, v(N+1; \lambda)$ is known. Of course, deg$\,F = 2N-2$. Given $F(\lambda)$ there are finitely
many possibilities for its $N$-degree factor $v(N+1; \lambda)$. However, in general, $v(N+1; \lambda)$ cannot be recovered uniquely from $F(\lambda)$
and one might suspect that, in general, it may not be possible to uniquely recover $L$ from $\{v(N; \lambda), u(N+1; \lambda)\}$. Actually, this
possibility can really happen, as it is demonstrated by the following (counter)example.

\medskip

\textbf{Example 1.} Let us take $N=4$ and consider the case $a(n) \equiv -1$ and $b(n)$ such that
\begin{equation}
b(1) = b(4) = \alpha + \sigma \sqrt{2},
\qquad
b(2) = b(3) = \alpha - \sigma \frac{\sqrt{2}}{2}
\label{H2}
\end{equation}
where $\alpha$ is any fixed real (or complex) number and $\sigma \in \{-1, 1\}$. Then \eqref{DB2} and \eqref{DB3} give
\begin{equation}
v(4; \lambda) = -u(5; \lambda) = -\lambda^3 + 3\alpha \lambda^2 - \left(3\alpha^2 - \frac{7}{2}\right) \lambda + \alpha^3 - \frac{7\alpha}{2}.
\label{H3}
\end{equation}
Hence, the sign $\sigma$ cannot be recovered from $v(4; \lambda)$ and $u(5; \lambda)$. In other words, there are two different potentials of period
$N=4$ corresponding to the same spectral data $\{v(4; \lambda), u(5; \lambda)\}$.

\medskip

Example 1 is, somehow, in contrast with \cite{L-G}.

We, now, wish to consider the following question: Suppose we are given a polynomial
\begin{equation}
D(\lambda) = (-1)^N \lambda^N + \sum_{k = 0}^{N-1} c_k \lambda^k.
\label{ND1}
\end{equation}
Is there an $N$-periodic operator $L_{\text{Schr}}$ whose discrete Hill discriminant is the given polynomial $D(\lambda)$?

Let us first give a lemma of algebraic flavor.

\medskip

\textbf{Lemma 1.} For $k = 1, \dots, N$ let $S_k(x_1, \dots, x_N)$ be the elementary symmetric polynomial in
the variables $x_1, \dots, x_N$ of degree $k$. Also, let $p_k(x_1, \dots, x_N)$, $k = 1, \dots, N$, be $N$
given polynomials in $x_1, \dots, x_N$ such that deg$\,p_k \leq k-1$. Then, the cardinality of the set $\Lambda$ of the distinct solutions
$(x_1, \dots, x_N) \in \mathbb{C}^N$ of the system of $N$ equations
\begin{equation}
S_k(x_1, \dots, x_N) = p_k(x_1, \dots, x_N),
\qquad
k = 1, \dots, N,
\label{L1}
\end{equation}
satisfies $1 \leq \#(\Lambda) \leq N!$.

\smallskip

\textit{Proof}. The result follows from the very simple observation that the system \eqref{L1} does not have solutions at infinity.

Let us make the above statement precise.
Consider the complex projective space $\mathbb{C}\mathbb{P}^N$  of (complex) dimension $N$. Recall that
$\mathbb{C}\mathbb{P}^N$ is the projective compactification of $\mathbb{C}^N$, which is constructed as follows:
The points of $\mathbb{C}\mathbb{P}^N$ have homogeneous
coordinates $(\xi_1, \dots, \xi_N, \xi_{N+1}) \ne (0, \dots, 0, 0)$ so that the set of coordinates
$(\lambda \xi_1, \dots, \lambda \xi_N, \lambda\xi_{N+1})$, $\lambda \in \mathbb{C} \setminus \{0\}$, represent
the same point of $\mathbb{C}\mathbb{P}^N$. The point of $\mathbb{C}\mathbb{P}^N$ with homogeneous coordinates
$(\xi_1, \dots, \xi_N, \xi_{N+1})$ with $\xi_{N+1} \ne 0$ can be identified with the point of $\mathbb{C}^N$ with
coordinates $(\xi_1 / \xi_{N+1}, \dots, \xi_N/ \xi_{N+1})$,  while the points of $\mathbb{C}\mathbb{P}^N$ with
homogeneous coordinates $(\xi_1, \dots, \xi_N, 0)$ are the so-called ``points at
infinity" and they do not correspond to points of $\mathbb{C}^N$.

We can, then, consider the ``projectified" \eqref{L1}, namely the system \eqref{L1} in $\mathbb{C}\mathbb{P}^N$, written in homogeneous coordinates:
\begin{equation}
S_k(\xi_1, \dots, \xi_N) = p_k(\xi_1 /\xi_{N+1}, \dots, \xi_N / \xi_{N+1}) \xi_{N+1}^k,
\qquad
k = 1, \dots, N,
\label{L2}
\end{equation}
so that the solutions of \eqref{L2} with $\xi_{N+1} \ne 0$ correspond to the solutions of \eqref{L1} in $\mathbb{C}^N$,
while the solutions of \eqref{L2} with $\xi_{N+1} = 0$ are the solutions at infinity.
Since for each $k = 1, \dots, N$ the polynomial $p_k(x_1, \dots, x_N)$ has degree $\leq k-1$, we
must have that all the quantities $p_k(\xi_1 /\xi_{N+1}, \dots, \xi_N / \xi_{N+1}) \xi_{N+1}^k$,
$k = 1, \dots, N$, vanish, if $\xi_{N+1} = 0$. Hence, if we set $\xi_{N+1} = 0$ in \eqref{L2} we get
\begin{equation}
S_k(\xi_1, \dots, \xi_N) = 0,
\qquad
k = 1, \dots, N,
\label{L3}
\end{equation}
and \eqref{L3} implies easily that $\xi_1 = \cdots = \xi_N = 0$, which is impossible by the definition of
homogeneous coordinates. Therefore, \eqref{L2} cannot have solutions at infinity and, consequently, $\Lambda$ is a compact subset of $\mathbb{C}^N$.
But, then, by the Noether's Normalization Theorem \cite{E} we can conclude that $\Lambda$ must be a finite set. Thus, the proof is finished by
invoking B\'{e}zout's Theorem \cite{Sh}.
\hfill $\blacksquare$

\medskip

We are now ready for the main result.

\medskip

\textbf{Theorem 1.} Let $c_0, \ldots, c_{N-1}$ be given complex numbers. Then, there exist at least one and at most $N!$ different
$N$-periodic potentials $b(n)$ for which the discrete Hill discriminant of the corresponding operator $L_{\text{Schr}}$ (see \eqref{L0}) is
\begin{equation}
\Delta(\lambda) = (-1)^N \lambda^N + \sum_{k = 0}^{N-1} c_k \lambda^k.
\label{ND2}
\end{equation}

\smallskip

\textit{Proof}. First, for typographical convenience let us write $b_n$ instead of $b(n)$ and, also, introduce the
notation
\begin{equation}
Q(n) := b_n - \lambda.
\label{E2}
\end{equation}
Next, let $\chi(n) = \chi(n; \lambda)$ and $\gamma(n) = \gamma(n; \lambda)$ be the (unique) solutions of $(L_{\text{Schr}} \, w)(n) = \lambda w(n)$ satisfying the initial conditions \eqref{B2a}, namely $\chi(0) = 1$, $\chi(1) = 0$ and $\gamma(0) = 0$, $\gamma(1) = 1$ respectively. In view of
\eqref{B3a}, the Wronskian (or Casoratian) of $\chi$ and $\gamma$ is
\begin{equation}
W[\chi,\gamma] :=
\left|
  \begin{array}{cc}
    \chi(n) & \gamma(n) \\
    \chi(n+1) & \gamma(n+1) \\
  \end{array}
\right|
\equiv 1,
\label{E2a}
\end{equation}
and this implies that the discriminant of $L_{\text{Schr}}$ is (as we have already seen in \eqref{B7a})
\begin{equation}
\Delta(\lambda) = \chi(N; \lambda) + \gamma(N+1; \lambda),
\label{E2b}
\end{equation}
By using the notation introduced above we can write the difference equation satisfied by $\chi(n)$ and $\gamma(n)$ as
\begin{equation}
w(n+1) = Q(n) w(n) - w(n-1),
\qquad
n \in \mathbb{Z}.
\label{E3}
\end{equation}
It follows by easy induction that, for $n \geq 2$ the quantity $\chi(n; \lambda)$ is a polynomial in $\lambda$ of degree $n-2$
(since $\chi(2; \lambda) = -1$), while $\gamma(n; \lambda)$ is a polynomial in $\lambda$ of degree $n-1$ having the form
\begin{align}
\gamma(n) = Q(1) Q(2) \cdots Q(n-1) \ - &\sum_{(j_1, \dots, j_{n-3}) \,\in \,\mathcal{S}_{n-3}} Q(j_1) \cdots Q(j_{n-3})
\nonumber
\\
+ &\sum_{(j_1, \dots, j_{n-5}) \,\in \,\mathcal{S}_{n-5}} Q(j_1) \cdots Q(j_{n-5}) - \cdots,
\label{E4}
\end{align}
where the $\mathcal{S}_k$'s are certain sets of $k$-tuples of distinct integers between $1$ and $n$. For example,
\begin{equation*}
\gamma(4) = Q(1) Q(2) Q(3) - Q(1) - Q(3)
\end{equation*}
and
\begin{equation*}
\gamma(5) = Q(1) Q(2) Q(3) Q(4) - Q(1) Q(2) - Q(1) Q(4) - Q(3) Q(4) + 1.
\end{equation*}
By using \eqref{E4} in \eqref{E2b} (together with the fact that deg$_{\lambda} \, \chi(n; \lambda) = n-2$) and then in \eqref{B7} we obtain
\begin{align}
\Delta(\lambda)
= Q(1) Q(2) \cdots Q(N) \ - &\sum_{(j_1, \dots, j_{N-2}) \,\in \,\mathcal{T}_{N-2}} Q(j_1) \cdots Q(j_{N-2})
\nonumber
\\
+ &\sum_{(j_1, \dots, j_{N-4}) \,\in \,\mathcal{T}_{N-4}} Q(j_1) \cdots Q(j_{N-4}) - \cdots,
\label{E7}
\end{align}
where the $\mathcal{T}_k$'s are certain sets of $k$-tuples of distinct integers between $1$ and $N$.

Let us now assume that $\Delta(\lambda)$ is as in \eqref{ND2}. Then, in view of \eqref{E7} and \eqref{E2} we must have
\begin{equation}
c_{N-1} = (-1)^{N-1} S_1(b_1, \dots, b_N),
\label{E8b}
\end{equation}
\begin{equation}
c_{N-2} = (-1)^N S_2(b_1, \dots, b_N) + P_2(b_1, \dots, b_N),
\label{E8c}
\end{equation}
\begin{equation}
c_{N-3} = (-1)^{N-1} S_3(b_1, \dots, b_N) + P_3(b_1, \dots, b_N),
\label{E8d}
\end{equation}
\begin{equation*}
\vdots
\end{equation*}
\begin{equation}
c_0 = S_N(b_1, \dots, b_N) + P_N(b_1, \dots, b_N),
\label{E8f}
\end{equation}
where $S_k(b_1, \dots, b_N)$, $k = 1, \dots, N$, are the elementary symmetric polynomials in
the variables $b_1, \dots, b_N$ of degree $k$, while $P_k(b_1, \dots, b_N)$, $k = 2, \dots, N$, are
given polynomials in $b_1, \dots, b_N$ such that deg$\,P_k \leq k-1$ (e.g., \eqref{DB4} implies that $P_2(b_1, \dots, b_N) \equiv (-1)^{N-1} N$). Therefore, we can apply Lemma 1 to the system
\eqref{E8b}--\eqref{E8f} and conclude that it has at least one and at most $N!$ distinct solutions $(b_1, \dots, b_N)$.
\hfill $\blacksquare$

\medskip

In Section 5 we give some examples of operators $L_{\text{Schr}}$ with a given discriminant.

Finally, let us mention that one can follow the approach used for proving Theorem 1 in order to show existence of operators $L_{\text{Schr}}$ whose
certain associated polynomial quantity, say $u(N, \lambda)$ (see \eqref{B2}), is of a given form.

\section{Periodic Jacobi operators whose spectrum is a closed interval}
We begin with some observations regarding the multiplier $r(\lambda)$.

As we have already mentioned, since $\Delta(\lambda)^2 - 4$ is a polynomial of even degree, namely $2N$, it follows that $\infty$ is not a branch point of $\sqrt{\Delta(\lambda)^2 - 4}$. Hence, in view of \eqref{B14} we have that
$\infty$ is not a branch point of $r(\lambda)$. Moreover,
from \eqref{B10} and \eqref{DB4} we have that one of the branches of $r(\lambda)$ satisfies
\begin{equation}
r_j(\lambda) = (-1)^N \lambda^N \left[1 + O\left(\frac{1}{\lambda}\right)\right]
\qquad \text{as }\;
\lambda \to\infty,
\label{DC1}
\end{equation}
while, due to the first equation in \eqref{B9} the other branch satisfies
\begin{equation}
r_k(\lambda) = (-1)^N \lambda^{-N} \left[1 + O\left(\frac{1}{\lambda}\right)\right]
\qquad \text{as }\;
\lambda \to\infty,
\label{DC2}
\end{equation}
where, $\{j, k\} = \{1, 2\}$.

Suppose $r(\lambda)$ has no branch points in $\mathbb{C}$. Then, $\Delta(\lambda)^2 - 4$ is the square of a polynomial, which in view of
\eqref{B10} implies that both $r_1(\lambda)$ and $r_2(\lambda)$ are polynomials. However, this is impossible, e.g., due to \eqref{DC2}. It follows that $r(\lambda)$ must have at least one branch point in $\mathbb{C}$. But, if $r(\lambda)$ has exactly one branch point $\eta \in \mathbb{C}$, then from \eqref{B14} we get that $[\Delta(\lambda)^2 - 4] / (\lambda - \eta)$ is the square of a polynomial, which is impossible since its degree is odd, namely $2N-1$ (actually, by the same argument we can deduce that the number of branch points of $r(\lambda)$ cannot be odd). Therefore $r(\lambda)$
has at least two branch points in $\mathbb{C}$.

\smallskip

The following theorem characterizes the spectrum of $L$ in the case where $r(\lambda)$ has exactly two branch points.

\medskip

\textbf{Theorem 2.} Suppose that the multiplier $r(\lambda)$ associated to \eqref{B1} has exactly two branch points
$\eta, \theta \in \mathbb{C}$.
Then, $\eta$ and $\theta$ are periodic or antiperiodic eigenvalues of $L$ satisfying
\begin{equation}
\left(\frac{\eta - \theta}{4}\right)^N = \pm1
\label{T1a}
\end{equation}
(in particular $|\eta - \theta| = 4$) and the spectrum of $L$ is the line segment joining $\eta$ and $\theta$, namely
\begin{equation}
\sigma(L) = \{\lambda \in \mathbb{C} \, : \, \lambda = \eta + (\theta - \eta)\, t, \ \ 0 \leq t \leq 1\}.
\label{T1b}
\end{equation}

\smallskip

\textit{Proof}. Since $r(\lambda)$ has only two branch points $\eta, \theta$, it follows from \eqref{B14} that it must have the form
\begin{equation}
r(\lambda) = \frac{\Delta(\lambda) + Q_0(\lambda) \sqrt{(\lambda - \eta)(\lambda - \theta)}}{2},
\label{C1}
\end{equation}
where $Q_0(\lambda)$ is a polynomial of degree $N-1$. Equation \eqref{C1} can be also written as
\begin{equation}
r_1(\lambda), r_2(\lambda)
= \frac{\Delta(\lambda) \pm Q_0(\lambda) \sqrt{(\lambda - \eta)(\lambda - \theta)}}{2}.
\label{C1a}
\end{equation}
Formula \eqref{C1} suggests that it will be more convenient, instead of $\lambda$ to work with the spectral parameter $\zeta$ where
\begin{equation}
\zeta + \zeta^{-1} :=  \frac{4\lambda - 2(\eta + \theta)}{\eta - \theta},
\label{DC3}
\end{equation}
so that
\begin{equation}
(\lambda - \eta)(\lambda - \theta)
= \frac{(\eta - \theta)^2}{16} \left(\zeta - \zeta^{-1}\right)^2.
\label{DC4}
\end{equation}
Notice that formula \eqref{DC3} tells us that for every $\lambda \in \mathbb{C} \setminus \{\eta, \theta\}$ we get exactly two
values of $\zeta$, say $\zeta_1$ and $\zeta_2$ with $\zeta_1 \zeta_2 = 1$ (hence $\zeta_1, \zeta_2 \ne 0$), while for $\lambda = \eta$ we only
get $\zeta = 1$ and for $\lambda = \theta$ we only get $\zeta = -1$. Conversely, for every $\zeta \in \mathbb{C} \setminus \{0\}$ we get exactly
one $\lambda$.

Using \eqref{DC4} in \eqref{C1a} yields
\begin{equation}
r_1(\lambda), r_2(\lambda)
= \frac{\Delta(\lambda) \pm Q_1(\lambda) \left(\zeta - \zeta^{-1}\right)}{2},
\quad \text{where }\;
Q_1(\lambda) := \frac{\eta - \theta}{4} \, Q_0(\lambda).
\label{DC5}
\end{equation}
Formula \eqref{DC5} implies that $r_1(\lambda) = R_1(\zeta)$ and $r_2(\lambda) = R_2(\zeta)$, where $R_j(\zeta)$,
$j = 1, 2$, are rational functions. On the other hand, by recalling the first equality in \eqref{B9}, namely
$r_1(\lambda) r_2(\lambda) \, \equiv 1$, we know that the multipliers are different from $0$ for every
$\lambda \in \mathbb{C}$, while, e.g., from \eqref{B10} we know that they cannot have poles in $\mathbb{C}$. Since, in view
of \eqref{DC3}, $\zeta = 0$ is the only value of $\zeta$ which does not correspond to a complex number $\lambda$, we
can conclude that $\zeta = 0$ is the only possible zero or pole
of $R_j(\zeta)$, $j = 1, 2$, and, consequently, since $R_1(\zeta) R_2(\zeta)  \equiv 1$, we must have
\begin{equation}
R_1(\zeta), R_2(\zeta) = c^{\pm 1} \zeta^{\pm d}
\qquad \text{or, equivalently,} \qquad
r_1(\lambda), r_2(\lambda) = c^{\pm 1} \zeta^{\pm d},
\label{DC6}
\end{equation}
for some complex constant $c \ne 0$ and some integer $d \geq 1$ (remember that $r_j(\lambda)$, $j = 1, 2$,
cannot be constant).

We continue by recalling that, since $\eta$ (and $\theta$) is a branch point of $r(\lambda)$, the equation
\eqref{B14} implies that $r(\eta) = \pm 1$ (and also $r(\theta) = \pm 1$). In particular, this tells us that $\eta$ and $\theta$
are periodic or antiperiodic eigenvalue of $L$. On the other hand, as we have already mentioned
if $\lambda = \eta$, then \eqref{DC3} gives that $\zeta = 1$ (while if $\lambda = \theta$,
then $\zeta = -1$). Thus, by using this particular value
of $\lambda$, namely $\lambda = \eta$, in \eqref{DC6} we obtain
\begin{equation}
c = \pm 1
\label{DC6a}
\end{equation}
and, consequently, $c^{\pm 1} = c$. Hence \eqref{DC6} becomes
\begin{equation}
r_1(\lambda), r_2(\lambda) = c \zeta^{\pm d}.
\label{DC6b}
\end{equation}

Next, we consider asymptotics as $\lambda \to \infty$. From \eqref{DC3} we have
\begin{equation}
\zeta \sim \frac{4}{\eta - \theta} \, \lambda
\qquad \text{or} \qquad
\zeta^{-1} \sim \frac{4}{\eta - \theta} \, \lambda.
\label{DC7}
\end{equation}
By using \eqref{DC7} in \eqref{DC6b} we get
\begin{equation}
r_1(\lambda), r_2(\lambda) \sim c \left(\frac{4}{\eta - \theta} \, \lambda\right)^{\pm d}
\qquad \text{as }\;
\lambda \to \infty.
\label{DC8}
\end{equation}
Then, by comparing \eqref{DC8} with \eqref{DC1} and \eqref{DC2} we obtain
\begin{equation}
d = N
\qquad \text{and} \qquad
\left(\frac{\eta - \theta}{4}\right)^N = (-1)^N c = \pm1.
\label{DC9}
\end{equation}

Finally, in view of \eqref{B15}, \eqref{DC3}, and \eqref{DC6b} we have that $\sigma(L)$ is as in \eqref{T1b}.
\hfill $\blacksquare$

\medskip

\textbf{Remark 5.} (i) In the case where $N$ is an odd integer, formula \eqref{T1a} can be written, without loss of generality (i.e. by interchanging the roles of $\eta$ and $\theta$) in the slightly simpler form
\begin{equation}
\left(\frac{\eta - \theta}{4}\right)^N = 1
\label{DC11}
\end{equation}

(ii) From \eqref{DC6b} (recall that $d=N$) and \eqref{DC3} we get
\begin{equation}
r(\lambda) = c \left[\frac{2 \lambda -(\eta + \theta) + 2 \sqrt{(\lambda - \eta)(\lambda - \theta)}}{\eta - \theta}\right]^N,
\label{DC12a}
\end{equation}verifying, e.g., that $r(\eta) = c$ and $r(\theta) = (-1)^N c$, where $c = \pm1$ as in the proof of Theorem 2. Therefore, if
$N$ is even, then $\eta$ and $\theta$ are both periodic or both antiperiodic eigenvalues of $L$, while if $N$ is odd, then one of the $\eta$, $\theta$ is a periodic eigenvalues of $L$, while the other is an antiperiodic eigenvalues. Furthermore, by using \eqref{DC12a} in
\eqref{B9} we get
\begin{equation}
c \Delta(\lambda)
= \left[\frac{2 \lambda -(\eta + \theta) + 2 \sqrt{(\lambda - \eta)(\lambda - \theta)}}{\eta - \theta}\right]^N
+ \left[\frac{2 \lambda -(\eta + \theta) - 2 \sqrt{(\lambda - \eta)(\lambda - \theta)}}{\eta - \theta}\right]^N.
\label{DC12}
\end{equation}
Since by \eqref{DC9} we have $(\eta - \theta)^N = (-1)^N 4^N c$, formula \eqref{DC12} can be also written as
\begin{equation}
(-1)^N \Delta(\lambda)
= \left[\frac{2 \lambda -(\eta + \theta) + 2 \sqrt{(\lambda - \eta)(\lambda - \theta)}}{4}\right]^N
+ \left[\frac{2 \lambda -(\eta + \theta) - 2 \sqrt{(\lambda - \eta)(\lambda - \theta)}}{4}\right]^N.
\label{DC12b}
\end{equation}
From \eqref{DC12b} we get
\begin{equation}
\Delta(\lambda)
= (-1)^N \left[\lambda^N - \frac{N (\eta + \theta)}{2} \lambda^{N-1}
+ \frac{(2N-3)N (\eta^2 + \theta^2) +2(2N-1)N \eta \theta}{16} \lambda^{N-2} + \cdots\right]
\label{DC12c}
\end{equation}
and by comparing \eqref{DC12c} with \eqref{DB4} we obtain (recalling \eqref{DA3})
\begin{equation}
B_0 = \frac{1}{N}\sum_{j=1}^N b(j) = \frac{\eta + \theta}{2}
\label{DC13a}
\end{equation}
and
\begin{equation}
\sum_{1 \leq j < k \leq N} b(j) \, b(k) - \sum_{j=1}^N a(j)^2
= \frac{(2N-3)N (\eta^2 + \theta^2) +2(2N-1)N \eta \theta}{16},
\label{DC13b}
\end{equation}
where, of course, $\eta$ and $\theta$ must satisfy \eqref{T1a}.
\hfill $\diamondsuit$

\medskip

We continue with a converse of Theorem 2.

\medskip

\textbf{Theorem 3.} Suppose that the spectrum $\sigma(L)$ is a simple piecewise smooth arc in the complex plane joining two (distinct) numbers $\eta$ and $\theta$. Then $\eta$ and $\theta$ are the only branch points of the multiplier $r(\lambda)$. Consequently, due to
Theorem 1, $\eta$ and $\theta$ must satisfy \eqref{T1a} and $\sigma(L)$ must be the line segment joining them, as displayed in \eqref{T1b}.

\smallskip

\textit{Proof}. By the spectral characterization given in \eqref{B15} we have that
\begin{equation}
\lambda \in \sigma(L)
\qquad \text{if and only if} \qquad
|r_1(\lambda)| = 1,
\label{C3}
\end{equation}
since, in view of \eqref{B9} we have $|r(\lambda)| = 1$ if and only if $|r_1(\lambda)| = 1$.

Now, let $\lambda_{\star} \in \mathbb{C}$ be a branch point of $r(\lambda)$. Then, as we have already seen
$r(\lambda_{\star}) = \pm 1$ and
$\lambda_{\star} \in \sigma(L)$. Suppose $\lambda_{\star} \neq \eta, \theta$. Then, there is an
$\varepsilon > 0$ such that the intersection $\beta := \sigma(L) \cap D(\lambda_{\star}; \varepsilon)$, where $D(\lambda_{\star}; \varepsilon)$ is
the closed disk of radius $\varepsilon$ centered at $\lambda_{\star}$, is a subarc of $\sigma(L)$ whose endpoints lie on the boundary of
$D(\lambda_{\star}; \varepsilon)$. In other words $D(\lambda_{\star}; \varepsilon) \setminus \beta$ has two components (i.e. $\beta$ separates $D(\lambda_{\star}; \varepsilon)$ in two pieces).

The expansion of $r_1(\lambda)$ about $\lambda = \lambda_{\star}$ has the form
\begin{equation}
r_1(\lambda) = \pm 1 + \left(\lambda - \lambda_{\star}\right)^{1/2}
\sum_{m=0}^{\infty} c_m \left(\lambda - \lambda_{\star}\right)^m,
\label{C4}
\end{equation}
where $c_m$ cannot vanish for all $m \geq 0$. Formula \eqref{C4} implies that the image of the arc $\beta$ under the map $r_1(\lambda)$
must be a piecewise smooth curve forming a right angle at $r_1(\lambda_{\star})$, hence it
cannot be a subset of the unit circle, and this is in contradiction with \eqref{C3}. Therefore, it is impossible to have
$\lambda_{\star} \neq \eta, \theta$. Hence, the only possible branch points of $r(\lambda)$ are $\eta$ and $\theta$. Since
$r(\lambda)$ must have at least two branch points, it follows that they have to be precisely $\eta$ and $\theta$.
\hfill $\blacksquare$

\medskip

Let us mention that we believe that the continuous versions of Theorems 2 and 3 are also valid, and that the proofs should follow the same
philosophy as their discrete counterparts (some perhaps ``harmless" differences are that instead of polynomials we have entire functions of order
$\leq 1/2$, $\eta = \infty$ or $\theta = \infty$, etc).

\medskip

\textbf{Remark 6.} From Theorems 2 and 3 it follows that if the multiplier $r(\lambda)$ associated to \eqref{B1} has exactly two branch points
$\eta, \theta$ or if the spectrum $\sigma(L)$ of $L$ is a simple piecewise smooth arc in the complex plane joining two complex numbers
$\eta$ and $\theta$, then there is an $\varpi$ with $\varpi^{2N} = 1$ such that
\begin{equation}
\sigma\left(\varpi L - \eta - 2\right) = [-2, 2].
\label{C5}
\end{equation}
Also, if $\sigma(L) = [-2, 2]$, then by Theorem 2 the only branch points of $r(\lambda)$ are $\eta = -2$ and $\theta = 2$. Hence
\eqref{DC12a} together with the asymptotic formulas \eqref{DC1} and \eqref{DC2} yield that
\begin{equation}
r(\lambda) = \left[\frac{-\lambda + \sqrt{\lambda^2 - 4}}{2}\right]^N = \tilde{r}(\lambda),
\label{C6}
\end{equation}
where (recall \eqref{DC0c6}) $\tilde{r}(\lambda)$ is the Floquet multiplier of the unperturbed case. Consequently,
the Hill discriminant $\Delta(\lambda)$ of $L$ must be equal to the unperturbed discriminant, namely (in view of \eqref{DC0c4})
\begin{equation}
\Delta(\lambda) = \tilde{\Delta}_N(\lambda) = z^N + z^{-N}
= \left(\frac{-\lambda +\sqrt{\lambda^2 - 4}}{2}\right)^N + \left(\frac{-\lambda -\sqrt{\lambda^2 - 4}}{2}\right)^N.
\label{C7}
\end{equation}
\hfill $\diamondsuit$

\medskip

\textbf{Example 2.} (i) If $N = 2$, $a(n) = i(-1)^n$, and $b(n) = 2(-1)^n$, then $\sigma(L) = [-2, 2]$.

(ii) If $N = 4$, $a(n) = (1+i) \, i^n/\sqrt{2}$, and $b(n) = (-1)^n \sqrt{2}$, then, again, $\sigma(L) = [-2, 2]$.

\subsection{Examples of discrete Schr\"{o}dinger operators whose spectrum is the interval $[-2, 2]$}
Suppose $L$ becomes the discrete Schr\"{o}dinger operator $L_{\text{Schr}}$ of \eqref{L0} and its spectrum is the closed interval $[-2, 2]$. Then,
by Remark 6 we know that its discriminant is given by \eqref{C7}. We can, therefore, apply Theorem 1 and conclude that there is at least one
and at most $N!$ such operators. Each of these operators is determined by its potential $(b_1, \dots, b_N)$, which in turn is a solution of
the system \eqref{E8b}--\eqref{E8f} (see also \eqref{E7}), where $c_j = \tilde{c}_j$, $j = 0, \dots, N-1$ are the coefficients of
$\tilde{\Delta}_N(\lambda)$.

Recall that from the discussion presented in Subsection 2.1 it follows that if $(b_1, \dots, b_N)$ is a solution of the system
\eqref{E8b}--\eqref{E8f}, namely if $b(n) = b_n$, $n = 1, \dots, N$, is a potential whose spectrum is $[-2, 2]$, then the same is true for
its ``cyclic permutations" $(b_2, b_3, \dots, b_N, b_1)$, $(b_3, b_4, \dots, b_N, b_1, b_2)$, etc. Also, since the coefficients
$\tilde{c}_j$, $j = 0, \dots, N-1$, as well as the coefficients of the system
\eqref{E8b}--\eqref{E8f} are rational, it follows that if $(b_1, \dots, b_N)$ is a solution, then so is $(\bar{b_1}, \dots, \bar{b_N})$ (where the bar denotes complex conjugation) and, furthermore all $b_1, \dots, b_N$ are algebraic numbers.

\medskip

\textbf{Example 3.} (i) For $N = 2$ (so that $\tilde{\Delta}_2(\lambda) = \lambda^2 - 2$) and $N = 3$
(so that $\tilde{\Delta}_3(\lambda) = -\lambda^3 + 3\lambda$) it is easy to check that the only (complex) solution of \eqref{E8b}--\eqref{E8f} is
the zero solution, namely $b(n) \equiv 0$.

\smallskip

(ii) For $N = 4$ (so that $\tilde{\Delta}_4(\lambda) = \lambda^4 - 4\lambda^2 + 2$) the system \eqref{E8b}--\eqref{E8f} becomes
\begin{equation*}
b_1 + b_2 + b_3 + b_4 = 0,
\end{equation*}
\begin{equation*}
b_1 b_2 + b_1 b_3 + b_1 b_4 + b_2 b_3 + b_2 b_4 + b_3 b_4 = 0,
\end{equation*}
\begin{equation*}
b_1 b_2 b_3 + b_1 b_2 b_4 + b_1 b_3 b_4 + b_2 b_3 b_4 = 0,
\end{equation*}
\begin{equation*}
b_1 b_2 b_3 b_4 =  b_1 b_2 + b_2 b_3 + b_3 b_4 + b_4 b_1.
\end{equation*}
It follows  that $b_1$, $b_2$, $b_3$, and $b_4$ are the roots of the equation $x^4 + \alpha = 0$, where
\begin{equation}
\alpha =  b_1 b_2 + b_2 b_3 + b_3 b_4 + b_4 b_1.
\label{E9e}
\end{equation}
Writing $b_1, b_2, b_3, b_4 = \pm(1 \pm i) \, \alpha^{1/4} / \sqrt{2}$ and substituting in \eqref{E9e} yields $\alpha = 0$ or $\alpha = 4$. From the
value $\alpha = 0$ we only get the obvious solution $b(n) \equiv 0$, whereas the value $\alpha = 4$ yields a total of eight distinct solutions:
\begin{equation}
\left[
  \begin{array}{c}
    b_1 \\
    b_2 \\
    b_3 \\
    b_4 \\
  \end{array}
\right]
=
\left[
  \begin{array}{c}
    1+i \\
    1-i \\
    -1+i \\
    -1-i \\
  \end{array}
\right],
\left[
  \begin{array}{c}
    -1-i \\
    1+i \\
    1-i \\
    -1+i \\
  \end{array}
\right],
\left[
  \begin{array}{c}
    -1+i \\
    -1-i \\
    1+i \\
    1-i \\
  \end{array}
\right],
\left[
  \begin{array}{c}
     1-i \\
    -1+i \\
    -1-i \\
    1+i \\
  \end{array}
\right]
\label{E9e2}
\end{equation}
and
\begin{equation}
\left[
  \begin{array}{c}
    b_1 \\
    b_2 \\
    b_3 \\
    b_4 \\
  \end{array}
\right]
=
\left[
  \begin{array}{c}
    1-i \\
    1+i \\
    -1-i \\
    -1+i \\
  \end{array}
\right],
\left[
  \begin{array}{c}
    -1+i \\
    1-i \\
    1+i \\
    -1-i \\
  \end{array}
\right],
\left[
  \begin{array}{c}
    -1-i \\
    -1+i \\
    1-i \\
    1+i \\
  \end{array}
\right],
\left[
  \begin{array}{c}
     1+i \\
    -1-i \\
    -1+i \\
    1-i \\
  \end{array}
\right].
\label{E9e3}
\end{equation}
Notice that the last three solutions in \eqref{E9e2} are the cyclic permutations of the first solution, while the four solutions
in \eqref{E9e3} are the complex conjugates of the solutions in \eqref{E9e2}. The first solution in \eqref{E9e2} corresponds to the potential
\begin{equation*}
b(n) =  -\frac{1+i}{2} \, i^n -i (-1)^n -\frac{1-i}{2} \, (-i)^n,
\end{equation*}
while (see Subsection 2.2) the other seven solutions correspond to the shifts of this potential, namely $b^1(n)$, $b^2(n)$, and $b^3(n)$, and to the
complex conjugates of those four potentials (changing $b(n)$ to $b^{\sharp}(n) = b(-n)$ does not produce any new solutions). All these eight
potentials, as well as the trivial potential $b(n) \equiv 0$ have spectrum $[-2, 2]$. Thus, there are only nine distinct solutions, while $4! = 24$.

\smallskip

(iii) For $N = 5$ (so that $\tilde{\Delta}_5(\lambda) = -\lambda^5 + 5\lambda^3 - 5\lambda$) the system \eqref{E8b}--\eqref{E8f} becomes
\begin{equation*}
S_1(b_1, b_2, b_3, b_4, b_5) = 0,
\end{equation*}
\begin{equation*}
S_2(b_1, b_2, b_3, b_4, b_5) = 0,
\end{equation*}
\begin{equation*}
S_3(b_1, b_2, b_3, b_4, b_5) = 0,
\end{equation*}
\begin{equation*}
S_4(b_1, b_2, b_3, b_4, b_5) = b_1 b_2 + b_2 b_3 + b_3 b_4 + b_4 b_5 + b_5 b_1,
\end{equation*}
\begin{equation*}
S_5(b_1, b_2, b_3, b_4, b_5) = b_1 b_2 b_3 + b_2 b_3 b_4 + b_3 b_4 b_5 + b_4 b_5 b_1 + b_5 b_1 b_2.
\end{equation*}
We can find some (nontrivial) solutions by looking for solutions such that $b_j = 0$ for some $j$, say $b_5 = 0$. Then, the system becomes
\begin{equation*}
S_1(b_1, b_2, b_3, b_4) = S_2(b_1, b_2, b_3, b_4) = S_3(b_1, b_2, b_3, b_4) = 0,
\end{equation*}
\begin{equation*}
b_1 b_2 b_3 b_4 = b_1 b_2 + b_2 b_3 + b_3 b_4,
\end{equation*}
\begin{equation*}
b_1 b_2 b_3 + b_2 b_3 b_4 = 0.
\end{equation*}
(although we have five equations with four unknowns, as we will see the resulting system has nine distinct solutions). If $b_2 b_3 = 0$, then we must have $b_j = 0$ for all $j = 1, \dots, 5$. Thus, let us assume $b_2 b_3 \ne 0$. In this case the last equation of the system can be simplified as
\begin{equation*}
b_1 + b_4 = 0.
\end{equation*}
As in the case (ii) it follows  that $b_1$, $b_2$, $b_3$, and $b_4$ are the roots of the equation $x^4 + \alpha = 0$, where
\begin{equation}
\alpha =  b_1 b_2 + b_2 b_3 + b_3 b_4.
\label{E11f}
\end{equation}
Writing $b_1, b_2, b_3, b_4 = \pm(1 \pm i) \, \alpha^{1/4} / \sqrt{2}$ and substituting in \eqref{E11f} yields $\alpha = 0$, $\alpha = 3+4i$,
or $\alpha = 3-4i$. From the value $\alpha = 0$ we only get the obvious solution $b(n) \equiv 0$.
The value $\alpha = 3+4i$ yields the solutions:
\begin{equation}
\left[
  \begin{array}{c}
    b_1 \\
    b_2 \\
    b_3 \\
    b_4 \\
    b_5 \\
  \end{array}
\right]
=
\left[
  \begin{array}{c}
    (1+i) \rho \\
    (1-i) \rho  \\
    -(1-i) \rho  \\
    -(1+i) \rho  \\
     0 \\
  \end{array}
\right],
\left[
  \begin{array}{c}
     (1-i) \rho \\
    -(1+i) \rho  \\
    (1+i) \rho  \\
    -(1-i) \rho  \\
     0 \\
  \end{array}
\right],
\left[
  \begin{array}{c}
   -(1-i) \rho \\
    (1+i) \rho  \\
    -(1+i) \rho  \\
    (1-i) \rho  \\
     0 \\
  \end{array}
\right],
\left[
  \begin{array}{c}
      -(1+i) \rho \\
     -(1-i) \rho  \\
    (1-i) \rho  \\
    (1+i) \rho  \\
     0 \\
  \end{array}
\right],
\label{E11g}
\end{equation}
where
\begin{equation}
\rho := \frac{\sqrt{\sqrt{5} + 2}}{2} + \frac{\sqrt{\sqrt{5} - 2}}{2} \, i.
\label{E11h}
\end{equation}
From the value $\alpha = 3-4i$ we get another set of four solutions, which are the complex conjugates of the solutions given in \eqref{E11g}.

An amusing observation is that these eight solutions can be also expressed as
\begin{equation*}
b_1 = \pm\frac{1}{\sqrt{\phi}} \pm  i\sqrt{\phi},
\qquad
b_2 = \pm i b_1,
\qquad
b_3 = -b_2,
\qquad
b_4 = -b_1,
\qquad
b_5 = 0
\end{equation*}
(for all eight different choices of the plus$/$minus signs), where $\phi$ is the golden ratio, i.e.
\begin{equation*}
\phi = \frac{1 + \sqrt{5}}{2}.
\end{equation*}
Finally, the cyclic permutations of the solutions in \eqref{E11g} and their complex conjugates produce a set of thirty two new solutions. The
transformation $b(n) \to b(-n)$ does not yield any new solutions. Thus we have found a total of forty distinct solutions, plus the obvious
(trivial) solution $b(n) \equiv 0$. Here we do not claim that we have found all the distinct solutions (since $5! = 120$, it is possible that
more solutions exist).

\subsection{The general operator}
We now consider again the more general Jacobi operator $L$ as introduced in \eqref{A3}. We are interested in the case where the spectrum $\sigma(L)$
is a simple piecewise smooth arc in the complex plane joining the numbers $\eta$ and $\theta$. Then, by Theorem 2 and Remark 6 we can assume,
essentially without loss of generality that $\sigma(L) = [-2, 2]$.

The following theorem is the discrete analog of a result of V. Guillemin and A. Uribe \cite{G-U2}.

\medskip

\textbf{Theorem 4.} Suppose that $\sigma(L) = [-2, 2]$. Then, the eigenvalues of $L|_{\mathcal{P}_{2N}}$, where $\mathcal{P}_{2N}$ is the
vector space of $2N$-periodic sequences as introduced in \eqref{G1}, or,
equivalently, the eigenvalues of the matrix $L_{2N}$ of \eqref{G3}, are the numbers given in \eqref{DC0f}, namely
\begin{equation}
\lambda_k = \tilde{\lambda}_k = -2 \cos\left(\frac{\pi k}{N}\right),
\qquad
k = 0, 1, \ldots, N.
\label{G4}
\end{equation}
Furthermore, for $k = 1, \dots, N-1$, either there are two linearly independent eigenfunctions (in $\mathcal{P}_{2N}$) corresponding to the eigenvalue $\lambda_k$ or there is a two-dimensional generalized eigenspace (subspace of $\mathcal{P}_{2N}$) of $L|_{\mathcal{P}_{2N}}$ associated to
$\lambda_k$.

\smallskip

\textit{Proof}. From Remark 6 we have that if $\sigma(L) = [-2, 2]$, then $r(\lambda) = \tilde{r}(\lambda)$, where $\tilde{r}(\lambda)$ is the
Floquet multiplier of the unperturbed case given by \eqref{DC0c6}. In particular, $r(\lambda) = \pm 1$ (i.e. $r_1(\lambda) = r_2(\lambda) = \pm 1$)
if and only if $\lambda = \lambda_k$ for some $\lambda_k$ of \eqref{G4} and, hence, the Floquet solutions of \eqref{B1} which are either $N$-periodic or $N$-antiperiodic are the solutions $\phi(n; \lambda_k)$, $k = 0, 1, \dots, N$, with $\lambda_k$ as in \eqref{G4}.

Now, as we have already seen in Subsection 2.1, a $2N$-periodic solution of \eqref{B1} is necessarily an $N$-periodic or $N$-antiperiodic
Floquet solution and vice versa.
Therefore, the eigenfunctions of $L|_{\mathcal{P}_{2N}}$ (in $\mathcal{P}_{2N}$) are precisely the $N$-periodic and $N$-antiperiodic Floquet
solutions and the spectrum of the operator $L|_{\mathcal{P}_{2N}}$ (whose matrix is $L_{2N}$) is given by \eqref{G4}.

Next, let $\lambda = \lambda_k$ for some $k =1, \dots, N-1$. If there are two linearly independent Floquet solutions $\phi_1(n; \lambda_k)$ and
$\phi_2(n; \lambda_k)$, then they are both $2N$-periodic and hence $L|_{\mathcal{P}_{2N}}$ has two linearly independent eigenfunctions
(in $\mathcal{P}_{2N}$).

Now, fix a $\lambda_k$, $k =1, \dots, N-1$, and suppose there is only one Floquet solution, say $\phi(n; \lambda_k)$ corresponding to $\lambda_k$,
normalized so that
\begin{equation}
\phi(0; \lambda_k) = 1
\label{G5}
\end{equation}
(formula \eqref{B14c} shows how to construct $\phi(n; \lambda_k)$; if $\phi(0; \lambda_k) = 0$ formula \eqref{B14c} fails, but then, instead of
\eqref{G5} we can normalize $\phi(n; \lambda_k)$ so that $\phi(1; \lambda_k) = 1$). As we have seen
(recall \eqref{B13}), in this case there is a solution $g(n)$ of \eqref{B1}, linearly independent of $\phi(n; \lambda_k)$, satisfying
\begin{equation}
(Sg)(n) = g(n + N) = r(\lambda_k) g(n) + \phi(n; \lambda_k)
\qquad \text{for all }\; n \in \mathbb{Z},
\label{G6}
\end{equation}
where $r(\lambda_k) = \pm1$.

Let us point out that, although $\lambda_k$ is not a branch point of $r(\lambda)$, the Floquet matrix $S(\lambda)$ has a Jordan anomaly at
$\lambda = \lambda_k$; in particular $\phi_1(n; \lambda_k) = \phi_2(n; \lambda_k)$. We say that $S(\lambda)$
\textit{has a pathology of the second kind over} $\lambda_k$ (this terminology was introduced in \cite{Pa}).

To continue we consider the system
\begin{equation}
L\phi(n; \lambda) = \lambda \phi(n; \lambda),
\quad
\phi(n + N; \lambda) = r(\lambda) \phi(n; \lambda),
\qquad
\lambda \in \mathbb{C},
\label{G7}
\end{equation}
where $\phi(n; \lambda)$ is as in \eqref{B14c}. Differentiating \eqref{G7} with respect to $\lambda$ yields
\begin{equation}
L\phi_{\lambda}(n; \lambda) = \lambda \phi_{\lambda}(n; \lambda) + \phi(n; \lambda),
\quad
\phi_{\lambda}(n + N; \lambda) = r(\lambda) \phi_{\lambda}(n; \lambda) + r'(\lambda) \phi(n; \lambda),
\label{G8}
\end{equation}
where the subscript $\lambda$ indicates derivative with respect to $\lambda$. If we fix $\lambda = \lambda_k$ and set (for typographical
convenience)
\begin{equation}
h(n) := \phi_{1 \lambda}(n; \lambda_k),
\label{G9}
\end{equation}
then \eqref{G8} becomes
\begin{equation}
(Lh)(n) = \lambda_k h(n) + \phi(n; \lambda_k),
\quad
h(n + N) = r(\lambda_k) h(n) + r_1'(\lambda_k) \phi(n; \lambda_k),
\label{G10}
\end{equation}
where $r_1'(\lambda_k) = \tilde{r}'(\lambda) = \pm (Ni/2) \sin(\pi k/N)$ and the sign depends on the branch of $r(\lambda_k)$; in fact,
$r_2'(\lambda_k) = -r_1'(\lambda_k)$). Finally, if we set
\begin{equation}
w(n) := h(n) - r_1'(\lambda_k) g(n),
\label{G11}
\end{equation}
then, in view of \eqref{G6} and \eqref{G10} (recall that $Lg = \lambda_k g$) we obtain
\begin{equation}
(Lw)(n) = \lambda_k w(n) + \phi(n; \lambda_k),
\qquad
w(n + N) = r(\lambda_k) w(n) = \pm w(n).
\label{G12}
\end{equation}
It follows that $w(n) \in \mathcal{P}_{2N}$ and $(L|_{\mathcal{P}_{2N}} - \lambda_k)^2 w(n) = 0$. Hence $\{\phi(n; \lambda_k), w(n)\}$ is a
generalized eigenspace of $L|_{\mathcal{P}_{2N}}$ associated to $\lambda_k$. Since the total dimension of the eigenspaces (pure or generalized),
for $k = 1, \dots, N-1$, is $2 (N-1) = 2N-2$ and we have two more eigenvalues of $L|_{\mathcal{P}_{2N}}$, namely $\lambda = \eta = -2$ and $\lambda = \lambda_N = 2$, we already have covered the $2N$-dimensional space $\mathcal{P}_{2N}$. Hence, the eigenspaces of
$\lambda = \eta = -2$ and $\lambda = \lambda_N = 2$ are one-dimensional, while to each $\lambda_k$, $k = 1, \dots, N-1$, corresponds a
two dimensional (pure or generalized) eigenspace.
\hfill $\blacksquare$

\medskip

\textbf{Remark 7.} A side product of Theorem 4 is that if $\sigma(L) = [-2, 2]$, then for $\lambda = -2$ and $\lambda = 2$ there is only one Floquet solution and consequently the Floquet matrices $S(-2)$ and $S(2)$ have a Jordan anomaly (at the same time, $r(\lambda)$ has a branch point at
$\lambda = \pm 2$; recall, also, that $r(-2) = 1$ and $r(2) = (-1)^N$).

\hfill $\diamondsuit$

\medskip

\textbf{Example 4.} Regarding the case $N = 4$: In the unperturbed case the matrix $L_8$ of \eqref{G3} is similar to the diagonal matrix
diag$\,[-2, -\sqrt{2}, -\sqrt{2}, 0, 0, \sqrt{2}, \sqrt{2}, 2]$. As
for the eight cases presented in the formulas \eqref{E9e2} and \eqref{E9e3} of Example 3(ii), the associated matrix $L_8$ is similar to the Jordan
canonical matrix
\begin{equation}
\left[
  \begin{array}{cccccccc}
    \mathbf{-2} & 0 & 0 & 0 & 0 & 0 & 0 & 0\\
    0 & \mathbf{-\sqrt{2}} & \mathbf{1} & 0 & 0 & 0 & 0 & 0\\
		0 & \mathbf{0} & \mathbf{-\sqrt{2}} & 0 & 0 & 0 & 0 & 0\\
		0 & 0 & 0 & \mathbf{0} & \mathbf{1} & 0 & 0 & 0\\
		0 & 0 & 0 & \mathbf{0} & \mathbf{0} & 0 & 0 & 0\\
		0 & 0 & 0 & 0 & 0 & \mathbf{\sqrt{2}} & \mathbf{1} & 0\\
		0 & 0 & 0 & 0 & 0 & \mathbf{0} & \mathbf{\sqrt{2}} & 0\\
        0 & 0 & 0 & 0 & 0 & 0 & 0 & \mathbf{2}\\
  \end{array}
\right].
\label{G13}
\end{equation}

\medskip

Finally, we present a Borg-type theorem for the general operator $L$ with complex coefficients.

\medskip

\textbf{Theorem 5.} Suppose that $\sigma(L) = [-2, 2]$ and that the matrix $L_{2N}$ of \eqref{G3} is diagonalizable (i.e. it has $2N$ linearly
independent pure eigenvectors). Then:

(i) If $N$ is odd, we must have $b(n) \equiv 0$ and $a(n)^2 \equiv 1$, i.e. $L$ is an essentially unperturbed operator (recall Definition 2).

(ii) If $N$ is even, say $N = 2M$, then $b(n) \equiv 0$ and $a(n)^2 = 1 + (-1)^n s$, where
$s^2 = 1 - e^{2k\pi i/M}$ for some $k \in \{0, 1, \ldots, M-1\}$.

\smallskip

\textit{Proof}. Setting $\eta = -2$ and $\theta = 2$ in formulas \eqref{DC13a} and \eqref{DC13b} of Remark 5 yields
\begin{equation}
N B_0 = \sum_{j=1}^N b(j) = 0
\label{BT1}
\end{equation}
and
\begin{equation}
\sum_{1 \leq j < k \leq N} b(j) \, b(k) - \sum_{j=1}^N a(j)^2 = -N,
\label{BT2}
\end{equation}
respectively. Furthermore, from Remark 6 we know that the assumption $\sigma(L) = [-2, 2]$ implies that the Hill discriminant of $L$ is
$\tilde{\Delta}_N(\lambda)$ of \eqref{DC0c4} and, consequently, the periodic$/$antiperiodic eigenvalues of $L$ are $\lambda_k = \tilde{\lambda}_k$
given by \eqref{DC0f}.

Since $L_{2N}$ is diagonalizable, Theorem 4 implies that to each $\lambda_k = \tilde{\lambda}_k$, $k = 1, \ldots, N-1$, correspond two linearly
independent Floquet solutions (i.e. we have coexistence), say $\phi_1(n; \lambda_k)$ and $\phi_2(n; \lambda_k)$, which are both periodic or
both antiperiodic. It follows that for each $\lambda_k$, $k = 1, \ldots, N-1$, there is a (nontrivial) linear combination
\begin{equation}
\phi(n; \lambda_k) = c_1(\lambda_k) \phi_2(n; \lambda_k) + c_2(\lambda_k) \phi_2(n; \lambda_k)
\label{H7}
\end{equation}
such that (since $\phi(n; \lambda_k)$ is either periodic or antiperiodic),
\begin{equation}
\phi(0; \lambda_k) = \phi(N; \lambda_k) = 0.
\label{H8}
\end{equation}
Therefore, the (distinct) numbers $\tilde{\lambda}_k$, $k = 1, \ldots, N-1$, are Dirichlet eigenvalues of $L$ and, since there are at most $N-1$
distinct Dirichlet eigenvalues (being the zeros of $\phi(b; \lambda)$) we must have that the Dirichlet spectrum of $L$ is
$\{\tilde{\lambda}_1, \dots, \tilde{\lambda}_{N-1}\}$. Thus,
the trace formula \eqref{D4b} becomes
\begin{equation}
\tilde{\lambda}_1 + \cdots + \tilde{\lambda}_{N-1} = B_0 N - b(0).
\label{H9}
\end{equation}
However, from \eqref{BT1} we know that $B_0 = 0$, while it is easy to check (e.g., from \eqref{DC0f}) that
$\tilde{\lambda}_1 + \cdots + \tilde{\lambda}_{N-1} = 0$.
Hence \eqref{H9} yields
\begin{equation}
b(0) = 0.
\label{H10}
\end{equation}
We now look at the ``shifted" operator $L^l$, where $l \in \{1, \dots, N\}$, and its associated operator $L^l$ (of course, $L^N = L$). As we have
discussed in Subsection 2.2 the Floquet solutions corresponding to $L^l$ are exactly the shifted Floquet solutions corresponding to $L$, while $L^l$
and $L$ have the same Floquet multiplier $r(\lambda)$. Thus $\sigma(L^l) = \sigma(L^l) = [-2, 2]$ and furthermore
$L^l_{2N}$ is diagonalizable. Therefore \eqref{H10} holds for $L^l$, namely
\begin{equation}
b(l) = b^l(0) = 0,
\qquad\text{for all }\; l = 1, \dots, N,
\label{H11}
\end{equation}
which means that $b(n) \equiv 0$.

We continue by noticing that if $b(n) \equiv 0$, then formula \eqref{BT2} becomes
\begin{equation}
\sum_{j=1}^N a(j)^2 = N.
\label{BT3}
\end{equation}
Now, the numbers $\tilde{\lambda}_1, \dots, \tilde{\lambda}_{N-1}$ (being the Dirichlet eigenvalues of $L$) are the zeros of the polynomial
$v(N; \lambda)$. Hence  formulas \eqref{DB3} and \eqref{DC0f} imply
\begin{equation}
-\sum_{j=1}^{N-2} a(j)^2 = \sum_{1 \leq j < k \leq N-1} \tilde{\lambda}_j \tilde{\lambda}_k = -(N-2).
\label{BT4}
\end{equation}
Thus, by using \eqref{BT4} in \eqref{BT3} we get
\begin{equation}
a(N-1)^2 + a(N)^2 = 2.
\label{BT5}
\end{equation}
In the very special case $N = 2$, formula \eqref{BT5} becomes $a(1)^2 + a(2)^2 = 2$, while the normalization \eqref{DA0} implies $a(1)^2 a(2)^2 = 1$.
Therefore, $a(1)^2 = a(2)^2 = 1$.

From now on we assume $N \geq 3$. Then, by considering the ``shifted" operator $L^l$ in place of $L$, $l = 1, \dots, N$, we can conclude from
\eqref{BT5} that
\begin{equation}
a(l)^2 + a(l+1)^2 = 2
\qquad\text{for all }\; l = 1, \dots, N
\label{BT6}
\end{equation}
(with $a(N+1) = a(1)$). Observe that \eqref{BT6} is a simple linear system of $N$ equations in $N$ unknowns, namely in $a(1)^2, \ldots, a(N)^2$.
By inspection, one solution of \eqref{BT6} is
\begin{equation}
a(1)^2 = a(2)^2 = \cdots = a(N)^2 = 1
\qquad \text{i.e. }\;
a(n)^2 \equiv 1.
\label{BT7}
\end{equation}
To find the other solutions of \eqref{BT6} (if there are any) we need to solve the associated homogeneous system
\begin{equation}
x_l + x_{l+1} = 0,
\quad
l = 1, \dots, N-1,
\qquad \text{and }\;
x_N + x_1 = 0.
\label{BT8}
\end{equation}
Suppose $x_1 = s$. Then $x_2 = -s$, $x_3 = s$, $\ldots$, $x_N = (1)^{N-1} s$ and, finally, $x_1 = (1)^N s$. Thus, if $N$ is odd, then we must have
$s = (1)^N s = -s$ and, consequently $s = 0$, which implies that the only solution of the homogeneous system \eqref{BT8} is the trivial solution
and, therefore, \eqref{BT6} implies that $a(n)^2 \equiv 1$.

It remains to examine the case of $N = 2M$. Here, the general solution of the  homogeneous system \eqref{BT8} is
$x_l = (1)^{l-1} s$, $l = 1, \ldots, N$, $s \in \mathbb{C}$. It follows that the general solution of \eqref{BT6} is
\begin{equation}
a(l)^2 = 1 + (1)^{l-1} s,
\qquad
l = 1, \dots, N.
\label{BT9}
\end{equation}
We, also, have the normalization condition \eqref{DA0} which implies
\begin{equation}
\prod_{l=1}^N a(l)^2 = 1.
\label{BT10}
\end{equation}
Substituting \eqref{BT9} in \eqref{BT10} yields (recall that $N = 2M$)
\begin{equation}
\left(1 - s^2\right)^M = 1,
\label{BT11}
\end{equation}
which tells us that $1 - s^2 = \rho$, equivalently $s^2 = 1 - \rho$, where $\rho$ is an $M$-th root of $1$.
\hfill $\blacksquare$

\medskip

Notice that in the case where $N = 2M$ the theorem implies that $a(n+2)^2 = a(n)^2$ for all $n \in \mathbb{Z}$.


Theorem 2 has a nice corollary.

\medskip

\textbf{Corollary 1.} If $a(n)$ and $b(n)$ are real-valued (equivalently, if $L$ is self-adjoint) and $\sigma(L) = [-2, 2]$, then $b(n) \equiv 0$
and $a(n)^2 \equiv 1$, i.e. $L$ is an essentially unperturbed operator.

\smallskip

\textit{Proof}. For real-valued $a(n)$ and $b(n)$ the matrix $L_{2N}$ of \eqref{G3} is real symmetric and hence diagonalizable. Therefore, the
corollary follows immediately from Theorem 5 since, even in the case $N = 2M$, the assumption that $a(n)$ is real forces $s$ to be $0$
(if $s^2$ is real, then $s^2 = 0$ or $s^2 = 2$; however, the latter cannot happen since it would make $a(n)^2$ strictly negative for certain values of $n$).
\hfill $\blacksquare$

\medskip

Corollary 1 is essentially not new (see \cite{F} or \cite{G-K-N-S}).

\section{Epilogue -- Some remarks on the Toda flow} Suppose that the coefficients  $a(n)$ and $b(n)$ of the operator $L$ depend on a parameter $t$ and that there is an
operator $B$ forming a Lax pair with $L$, namely $L$ and $B$ satisfy the equation
\begin{equation}
L_t = [B, L] := BL - LB,
\label{LP1}
\end{equation}
where, as usual, the subscript $t$ denotes derivative with respect to $t$. Then, as it is well known, the family of operators
$L = L(t)$, $t \in \mathbb{C}$, is isospectral, in the sense that the $l^2(\mathbb{Z})$-spectrum $\sigma(L)$ is independent of $t$. Actually,
the discriminant $\Delta(\lambda)$ of $L$ is independent of $t$.

One famous case of such an isospectral flow is the Toda flow, which is obtained by taking \cite{D-T}, \cite{G-H-T}
\begin{equation}
(B w)(n) := a(n) \, w(n+1) - a(n-1) \, w(n-1) = a_n \, w(n+1) - a_{n-1} \, w(n-1),
\quad
n \in \mathbb{Z}.
\label{LP2}
\end{equation}
In this case \eqref{LP1} can be written equivalently as \cite{D-T}
\begin{equation}
\left[a_n(t)^2\right]' = 2a_n(t)^2 \left[b_{n+1}(t) - b_n(t) \right],
\qquad
b_n'(t) = 2 \left[a_n(t)^2 - a_{n-1}(t)^2 \right],
\label{LP3}
\end{equation}
where the prime denotes derivative with respect to $t$.

Let us consider the case $\sigma(L) = [-2, 2]$, namely (recall Remark 6) $\Delta(\lambda) = \tilde{\Delta}_N(\lambda)$. Then, the Dirichlet
eigenvalues $\mu_1, \ldots, \mu_{N-1}$ of $L$ (see Subsection 2.3) satisfy the evolution equations \cite{D-T}
\begin{equation}
\mu_j' = -2 \sqrt{\tilde{\Delta}_N\left(\mu_j\right)^2 - 4} \prod_{k=1,\ k\ne j}^{N-1}\left(\mu_j - \mu_k \right)^{-1},
\qquad
j = 1, \ldots, N-1,
\label{LP4}
\end{equation}
where $\tilde{\Delta}_N(\lambda)$ is given by \eqref{DC0c4}.

One expects that there exist quantities similar to the "reflection coefficients" appearing in \cite{G-U1} for the continuous
periodic Schr\"{o}dinger operator $H$ with $\sigma(H) = [0, \infty)$, which determine $L$ and whose evolution under the Toda flow is very simple.

\bigskip

\textbf{Acknowledgments.} The author wishes to thank Professor Fritz Gesztesy for various helpful comments and suggestions, and Professor
Leonid Golinskii for pointing out certain relevant references.


\begin{thebibliography}{99}

%
%
%
%
%
%
%
%
%
%
%
%
%
%
%
%
%
%
%
%
%
%
%
%
%
%

\bibitem{B} Bender, C.M., Introduction to $\mathcal{PT}$-Symmetric Quantum Theory, arXiv:quant-ph/0501052v1 (2005), DOI: 10.1080/00107500072632


\bibitem{Bo} Borg, G., Eine Umkehrung der Sturm-Liouvilleschen Eigenwertaufgabe, \textit{Acta Math.}, \textbf{78}, 1--96 (1946)



\bibitem{D-T} Date, E. and Tanaka, S., Analog of Inverse Scattering Theory for the Discrete Hill’s Equation
and Exact Solution for the Periodic Toda Lattice, \textit{Progress of Theoretical Physics}, \textit{55}
(No. 2)  457--465 (February 1976)


\bibitem{E} Eisenbud, D., \textit{Commutative Algebra with a View Toward Algebraic Geometry},
Graduate Texts in Mathematics \textbf{150}, Springer-Verlag New York, 1995


\bibitem{F} Flaschka, H., Discrete and Periodic Illustrations of some Aspects of the Inverse Method, \textit{Lecture Notes in Physics}, \textbf{38},
441--466 (1975)

\bibitem{G} Gasymov, M.G., Spectral Analysis of a Class of Second-Order
Non-Self-Adjoint Differential Operators, \textit{Functional Anal. Appl.}, \textbf{14}, 11--15 (1980)

\bibitem{G-H} Gesztesy, F. and Holden, H., \textit{Soliton Equations and Their Algebro-Geometric Solutions},
\textit{Volume I: $(1 + 1)$-Dimensional Continuous Models}, Cambridge University Press, 2003

\bibitem{G-H-T} Gesztesy, F., Holden, H., and Teschl, G., The algebro-geometric Toda hierarchy initial value problem for complex-valued initial data,
\textit{Revista Matematica Iberoamericana}, \textbf{24} (no. 1), 117--182 (2008)

\bibitem{G-T} Gesztesy, F. and Tkachenko, V.A., When is a Non-Self-Adjoint Hill Operator a Spectral Operator of Scalar Type?,
\textit{C. R. Math. Acad. Sci. Paris}, \textbf{343} (4), 239--242 (2006),
http://dx.doi.org/10.1016/j.crma.2006.06.014

\bibitem{G-W} Gesztesy, F. and Weikard, R., Floquet theory revisited. Pages 67--84 of: Knowles, I. (ed), \textit{Differential Equations and
Mathematical Physics}, Boston: International Press, 1995

\bibitem{G-K-N-S} Golinskii, L., Kumar, K., Namboodiri, M.N.N., and Serra-Capizzano, S., A Note on a Discrete Version of Borg's Theorem via
Toeplitz-Laurent Operators with Matrix-Valued Symbols, \textit{Bollettino U.M.I.} (9) \textbf{VI}, 205--218 (2013)

\bibitem{G-U1} Guillemin, V. and Uribe, A., Hardy Functions and the Inverse Spectral Method, \textit{Comm. Partial Differential Equations},
\textbf{8}, 1455--1474 (1983)

\bibitem{G-U2} Guillemin, V. and Uribe, A., Spectral Properties of a Certain Class of Complex Potentials, \textit{Transactions of the American
Mathematical Society}, \textbf{279} (No. 2), 759--771 (1983)

\bibitem{Ho} Hochstadt, H., On the theory of Hill's matrices and related inverse spectral problems,
\textit{Linear Algebra and Appl.}, \textbf{11}, 41--52 (1975)

%
%

%




\bibitem{L-G} Levitan, B.M. and Gasymov, M.G., Determination of a differential equation by two of its spectra,
\textit{Russian Mathematical Surveys}, \textbf{19} (2), 1--63 (1964)

\bibitem{L} Last, Y., On the measure of gaps and spectra for discrete $1D$ Schr\"{o}dinger operators.
\textit{Commun. Math. Phys.}, \textbf{149}, 347--360 (1992)


%

\bibitem{N1} Na\u{\i}man, P.B., On the theory of periodic and limit-periodic Jacobian matrices,
\textit{Sov. Math. Dokl.}, \textbf{3}, 383--385 (1962)

\bibitem{N2} Na\u{\i}man, P.B., On the spectral theory of non-symmetric periodic Jacobi matrices,
\textit{Zap. Meh.-Mat. Fak. Har\c{p}rime kov. Gos. Univ. i Har\c{p}rime kov. Mat. Ob\v{s}\v{c}.}, \textbf{30} (4), 138--151 (1964). (Russian)

\bibitem{Pa} Papanicolaou, V.G., Some Results on Ordinary Differential Operators with Periodic Coefficients, \textit{Complex Anal. Oper. Theory},
\textbf{10}, 1227--1265 (2016) DOI 10.1007/s11785-015-0498-z

\bibitem{P-D} Papanicolaou, V.G. and Doumas, A.V., The Discrete Inverse Transmission Eigenvalue  Problem, {\it Inverse Problems} 27 015004 (2011)

\bibitem{P-T1} Pastur, L.A. and Tkachenko, V.A., Spectral Theory of Schr\"{o}dinger
Operators with Periodic Complex-Valued Potentials, \textit{Functional Anal. Appl.}, \textbf{22}, 156--158 (1988)

\bibitem{P-T2} Pastur, L.A. Tkachenko, V.A., Geometry of the Spectrum of the One-Dimensional Schr\"{o}dinger Equation with a Periodic Complex-Valued
Potential, \textit{Math. Notes}, \textbf{50}, 1045--1050 (1991)

\bibitem{P-T3} Pastur, L.A. and Tkachenko, V.A., An Inverse Problem for a Class of One-Dimensional Schr\"{o}dinger Operators with a Complex Periodic
Potential, \textit{Math. USSR-Izv.}, \textbf{37}, 611--629 (1991)

\bibitem{P} Peherstorfer, F., Orthogonal and extremal polynomials on several intervals,
\textit{Journal of Computational and Applied Mathematics}, \textbf{48}, 187--205 (1993)


\bibitem{RB} Rofe-Beketov, F.S., The Spectrum of Non-selfadjoint Differential Operators with Periodic Coefficients,
\textit{Soviet Math. Dokl.}, \textbf{4}, 1563--1566 (1963)

\bibitem{S-T} Sansuc, J.-J. and Tkachenko, V.A., Spectral Parametrization of Non-Selfadjoint Hill's Operators,
\textit{Journal of Differential Equations}, \textbf{125}, 366--384 (1996), Article no. 0035

\bibitem{S} Serov, M.I., Certain Properties of the Spectrum of a Non-selfadjoint Differential Operator of the Second Order,
\textit{Soviet Math. Dokl.}, \textbf{1}, 190--192 (1960)

\bibitem{Sh} Shafarevich, I.R., \textit{Basic Algebraic Geometry 1, Varieties in Projective Space}, Third Edition, Springer-Verlag Berlin Heidelberg,
2013

\bibitem{Si} Simon, B., \textit{Szeg˝\H{o}'s Theorem and Its Descendants: Spectral Theory for $L^2$ Perturbations of Orthogonal Polynomials},
Princeton University Press, Princeton and Oxford, 2011

\bibitem{T} Teschl, G., \textit{Jacobi Operators and Completely Integrable Nonlinear Lattices}, Mathematical Surveys and Monographs, Volume 72,
American Mathematical Society, 2000

\bibitem{T1} Tkachenko, V.A., Spectral Analysis of the One-Dimensional Schr\"{o}dinger Operator with a Periodic Complex-Valued Potential,
\textit{Soviet Math. Dokl.}, \textbf{5}, 413--415 (1964)

\bibitem{T2} Tkachenko, V.A., Spectral Analysis of a Nonselfadjoint Hill Operator, \textit{Dokl. Akad. Nauk SSSR}, \textbf{322} (No. 2) (1992);
\textit{Soviet Math. Dokl.}, \textbf{45} (No. 1), 78--82 (1992)

\bibitem{T3} Tkachenko, V.A., Discriminants and Generic Spectra of Nonselfadjoint Hill Operators, \textit{Adv. Sov. Math.}, \textbf{19},
41--71 (1994)

\bibitem{T4} Tkachenko, V.A., Spectra of non-selfadjoint Hill’s operators and a class of Riemann surfaces, Ann. Math. 143, 181--231 (1996)

\bibitem{V-TD} O. A. Veliev, O.A. and Toppamuk Duman, M., The Spectral Expansion for a Nonself-adjoint Hill Operator with a Locally Integrable
Potential, \textit{Journal of Mathematical Analysis and Applications}, \textbf{265}, 76--90 (2002)
doi:10.1006/jmaa.2001.7693, available online at http://www.idealibrary.com

\bibitem{Z} Zhernakov, N.V., Direct and inverse problems for a periodic Jacobi matrix,
\textit{Ukrainskii Matematicheskii Zhurnal}, \textbf{38} (No. 6), 785--788 (1986) (in Russian).
English translation: \textit{Ukrainian Math J.}, \textbf{38} (No. 6), 665--668 (1986)


\end{thebibliography}
\end{document}